\documentclass[12pt]{amsart}
\usepackage{amssymb,amsmath,amsthm,mathrsfs,multirow,xcolor,framed,url}
\def\bmaplesess{} % begin maple session
\def\emaplesess{} % end maple session
\def\bmaple{} % begin maple input
\def\emaple{} % end maple input
\def\bmapleout{} %begin maple output
\def\bmapleout{} %begin maple output
\def\emapleout{} %end maple output
\newcommand{\showtfile}[1]{}
\newtheorem{theorem}{Theorem}[section]
\newtheorem{lemma}[theorem]{Lemma}
\newtheorem{cor}[theorem]{Corollary}

\theoremstyle{definition}

\theoremstyle{remark}
\newtheorem{remark}[theorem]{Remark}

\numberwithin{equation}{section}
\newcommand{\abs}[1]{\left|{#1}\right|}

\usepackage{color}
\usepackage{tikz}
\usetikzlibrary{matrix}
\definecolor{shadecolor}{rgb}{0.8, 0.8, 0.8}
\newcommand\eqn[1]{(\ref{eq:#1})}
\def\maple{\textsc{maple \ }}
\def\maplep{\textsc{maple}.\ }
\def\sn{\smallskip\noindent}
\def\carrot{${}^\wedge$}
\newcommand\exampleshade[1]{
\begin{shaded}
\null\noindent
EXAMPLE:\\
{\tt
\noindent
#1
}
\end{shaded}
}
\newcommand\mapleshade[1]{
\begin{shaded}
\null\noindent
\hphantom{EXAMPLE:}\\
{\tt
\noindent
#1
}
\end{shaded}
}
\newcommand{\SL}{\mbox{SL}}

\newcommand\mylabel[1]{\label{#1}}
\newcommand{\beqs}{\begin{equation*}}
\newcommand{\eeqs}{\end{equation*}}
\newcommand{\beq}{\begin{equation}}
\newcommand{\eeq}{\end{equation}}
\newcommand\thm[1]{\ref{thm:#1}}

\newcommand\corol[1]{\ref{cor:#1}}

\newcommand\subsect[1]{\ref{subsec:#1}}

\newcommand{\C}{\mathbb{C}}

\newcommand{\Z}{\mathbb{Z}}
\newcommand{\Hup}{\mathcal{H}}

\renewcommand{\Im}{\mathrm{Im}\,}
\newcommand{\leg}[2]{\genfrac{(}{)}{}{}{#1}{#2}}

\DeclareMathOperator{\ord}{ord}
\DeclareMathOperator{\ORD}{ORD}
\DeclareMathOperator{\Up}{U}
\newcommand\Lpar[1]{\left(#1\right)}
\newcommand\umin[1]{\underset{#1}{\min}}
\begin{document}
%%%%%%%%%%%%%%%%%%%%%%%%%%%%%%%%%%%%%%%%%%%%%%%%%%%%%%%%%%%%%%%%%%%%%%%%%%%%%%%%%%%%%%%%%
\title{A tutorial for the MAPLE ETA package}                 

\author{Frank Garvan}
%%\address{Jie Frye,\\
%%Bunker Hill Community College,	Boston, MA 02129\\
%%\email{jiefrye@gmail.com}\\
%%Frank Garvan,\\
%%Department of Mathematics, University of Florida, Gainesville, FL 32601\\
%%\email{fgarvan@ufl.edu}
%%

%\thanks{ A preliminary version of this paper was presented by J.\ Frye
%on January 10, 2013 at JMM2013, San Diego. 
%F.\ Garvan was supported in part by a grant from
%the Simon's Foundation (\#318714). }


\thanks{The author was supported in part by a grant from 
the Simon's Foundation (\#318714).}
\date{July 9, 2019}

\maketitle
%%\date{November 22, 2013}
%%\date{May 22, 2018}
\date{July 9, 2019}

\begin{abstract}
This is a tutorial for using \texttt{ETA}, a MAPLE package for calculating
with Dedekind's eta function.
The \texttt{ETA} package is designed for proving
eta-product identities using the valence formula
for modular functions. 
\end{abstract}
%%%%%%%%%%%%%%%%%%%%%%%%%%%%%%%%%%%%%%%%%%%%%%%%%%%%%%%%%%%%%%%%%%%%%%%%%%%%%%%%%%%%%%%%%

%%\title*{Harmonic Sums\index{Harmonic Sums}, Polylogarithms\index{Polylogarithms},
%%Special Numbers\index{Special Numbers}, and their
%%Generalizations}
%%% Use \titlerunning{Short Title} for an abbreviated version of
%%% your contribution title if the original one is too long
%%\author{Jakob Ablinger and Johannes Bl\"umlein}
%%% Use \authorrunning{Short Title} for an abbreviated version of
%%% your contribution title if the original one is too long   
%%\institute{Jakob Ablinger,\\  Research Institute for Symbolic Computation (RISC)
%%Johannes Kepler University, Altenbergerstra\ss{}e 69, A-4040 Linz, Austria,
%%\email{jablinge@risc.uni-linz.ac.at},\\
%%Johannes Bl\"umlein, \\ Deutsches Elektronen-Synchrotron, DESY, Platanenallee 6, D-14738
%%Zeuthen, Germany,
%%\email{Johannes.Bluemlein@desy.de}}
%%%
%%% Use the package "url.sty" to avoid
%%% problems with special characters
%%% used in your e-mail or web address
%%%
%%\maketitle

%
%  --> relevant terms, to appear in the index shall alwys be put into \index{....}
%
\renewcommand{\theequation}{\arabic{section}.\arabic{equation}}
\allowdisplaybreaks

%%SECTION 1 %%%%%%%%%%%%%%%%%%%%%%%%%%%%%%%%%%%%%%%%%%%%%%%%%%%%%%%%%%%%%%%%%%%%%%%
\section{Introduction}
\setcounter{equation}{0}
\mylabel{sec:intro}

The Dedekind eta-function is defined by
\[ 
\eta(\tau) = q^{\frac{1}{24}} \prod_{n=1}^\infty (1-q^n),
\]
where $\tau \in \Hup:= \{ \tau \in \C : \Im \tau > 0\}$ 
and $q := e^{2\pi i \tau}$.

The main goal of the \texttt{ETA} \maple package is
to 
automatically prove identities for eta-products.

%%%%%%%%%%%%%%%%%%%%%%%%%%%%%%%%%%%%%%
\subsection{Installation Instructions}
\mylabel{subsec:install}

    First install the \texttt{qseries} package from 
\begin{center}
\url{http://qseries.org/fgarvan/qmaple/qseries}
\end{center}
 and  follow the directions on that page. Before proceeding it is advisable
to become familiar with the functions in the \texttt{qseries} package.
See \cite{Ga99b} for a tutorial.
Then go to
\begin{center}
\url{http://qseries.org/fgarvan/qmaple/ETA}
\end{center}
 to install the \texttt{ETA}
package.

%%SECTION 2 %%%%%%%%%%%%%%%%%%%%%%%%%%%%%%%%%%%%%%%%%%%%%%%%%%%%%%%%%%%%%%%%%%%%%%%
\section{Modularity of eta-products}
\mylabel{sec:modeta}
\setcounter{equation}{0}

          To prove a given eta-function identity one needs to  basically
do the following.
\begin{enumerate}
\item[(i)]
Rewrite the identity in terms of generalized eta-functions.
\item[(ii)]
Check that each term in the identity is a modular function on some
group $\Gamma_0(N)$.
\item[(iii)]
Determine the order at each cusp of $\Gamma_0(N)$ of each term
in the identity.
\item[(iv)]
Use the valence formula to determine up to which power of $q$ is needed to verify
the identity.
\item[(v)]
Finally prove the identity by carrying out the verification.
\end{enumerate}

In this section we explain how to carry out each of these steps in \maplep
Then we show how the whole process of proof can be automated.

%%%%%%%%%%%%%%%%%%%%%%%%%%%%%%%%%%%%
\subsection{Encoding eta-functions}     
\mylabel{subsec:encoding}

We encode $\eta(\tau)$ by \texttt{eta(tau)}. We will consider
eta-products of the form

\begin{equation}
f(\tau) = \prod_{d\mid N} \eta(d\tau)^{m_{d}},
\mylabel{eq:etapdefa}
\end{equation}
where $N$ is a positive integer and each $d>0$ and $m_{d}\in\Z$.
We encode the product
\begin{equation}
\prod_{j=1}^m \eta(t_j \tau)^{r_j}
\mylabel{eq:etatrdef}
\end{equation}
by the list
\begin{center}
[$t_1$, $r_1$, $t_2$, $r_2$, \dots $t_m$, $r_m$].
\end{center}
We call such a product an eta-product.
Such a list is called a generalized permutation and %KONDOref
usually written symbolically as
$$
t_1^{r_1} \, t_2^{r_2} \, \cdots \, t_m^{r_m}.
$$
We use the functions \texttt{GPmake} and \texttt{gp2etaprod} to convert
between symbolic forms of eta-products.

\noindent
\texttt{GPmake(etaprod)} --- converts an eta-product to a generalized
permutation.

\noindent
\texttt{gp2etaprod(gp)} --- converts the gneralized permutation \texttt{gp}
to an eta-product.

%%%%%%%%%%%%%%%%%%%%%%%%%%%%%%%%%%%%%%%%%%%%%%%%%%%%%%%%%
\showtfile{prog1.tex}
\exampleshade{
\sn
{\tt
\bmaplesess
\bmaple
$>$ \enskip   with(qseries): \newline
\emaple
\bmaple
$>$ \enskip   with(ETA): \newline
\emaple
\bmaple
$>$ \enskip   HM5:=q*mul((1-q\carrot(5*n))\carrot6/(1-q{\carrot}n)\carrot6,n=1..100): \newline
\emaple
\bmaple
$>$ \enskip   ep:=etamake(HM5,q,50);
\emaple
\bmapleout
$$
{\frac { \left( \eta \left( 5\,\tau \right)  \right) ^{6}}{ \left( \eta  \left( \tau \right)  \right) ^{6}}}
$$
\emapleout
\bmaple
$>$ \enskip   gp:=GPmake(ep);
\emaple
\bmapleout
$$
[5,6,1,-6]
$$
\emapleout
\bmaple
$>$ \enskip   gp2etaprod(gp);
\emaple
\bmapleout
$$
{\frac { \left( \eta \left( 5\,\tau \right)  \right) ^{6}}{ \left( \eta  \left( \tau \right)  \right) ^{6}}}
$$
\emapleout
\emaplesess
}\vskip0pt\noindent
}
%%%%%%%%%%%%%%%%%%%%%%%%%%%%%%%%%%%%%%%%%%%%%%%%%%%%%%%%%
%%     |\^/|     Maple 2018 (X86 64 WINDOWS)
%% ._|\|   |/|_. Copyright (c) Maplesoft, a division of Waterloo Maple Inc. 2018
%%  \  MAPLE  /  All rights reserved. Maple is a trademark of
%%  <____ ____>  Waterloo Maple Inc.
%%       |       Type ? for help.
%% >   with(qseries):
%% >   with(ETA):
%% >   HM5:=q*mul((1-q^(5*n))^6/(1-q^n)^6,n=1..100):
%% >   ep:=etamake(HM5,q,50);
%%                                                6
%%                                      eta(5 tau)
%%                                ep := -----------
%%                                               6
%%                                       eta(tau)
%% 
%% >   gp:=GPmake(ep);
%%                               gp := [5, 6, 1, -6]
%% 
%% >   gp2etaprod(gp);
%%                                             6
%%                                   eta(5 tau)
%%                                   -----------
%%                                            6
%%                                    eta(tau)
%% 
%% > quit
%% memory used=4.3MB, alloc=8.3MB, time=0.08
% 
% end prog1.tex
%%%%%%%%%%%%%%%%%%%%%%%%%%%%%%%%%%%%%%%%%%%%%%%%%%%%%%%%%

\subsection{$q$-Expansions}               
\mylabel{subsec:qexp}
There are two functions to compute $q$-expansions of eta-products.

\noindent
\texttt{etaprodtoqseries(etaprod)} --- returns the $q$-expansion of the 
eta-product \texttt{etaprod} up to $q^T$.

\noindent
\texttt{etaprodtoqseries2(etaprod)} --- a version of 
\texttt{etaprodtoqseries} that omits the factor $q^{1/24}$ in the
$q$-expansion of $\eta(\tau)$.

%%%%%%%%%%%%%%%%%%%%%%%%%%%%%%%%%%%%%%%%%%%%%%%%%%%%%%%%%
\showtfile{prog2.tex}
\exampleshade{
\sn
{\tt
\bmaplesess
\bmaple
$>$ \enskip  with(qseries): \newline
\emaple
\bmaple
$>$ \enskip  with(ETA): \newline
\emaple
\bmaple
$>$ \enskip  gp:=[2,2,1,-1]: \newline
\emaple
\bmaple
$>$ \enskip  ep:=gp2etaprod(gp);
\emaple
\bmapleout
$$
{\frac { \left( \eta \left( 2\,\tau \right)  \right) ^{2}}{\eta \left(  \tau \right) }} 
$$
\emapleout
\bmaple
$>$ \enskip  etaprodtoqseries(ep,50);
\emaple
\bmapleout
$$
\sqrt [8]{q}+{q}^{{\frac{9}{8}}}+{q}^{{\frac{25}{8}}}+{q}^{{\frac{49}{8 }}}+{q}^{{\frac{81}{8}}}+{q}^{{\frac{121}{8}}}+{q}^{{\frac{169}{8}}}+{q }^{{\frac{225}{8}}}+{q}^{{\frac{289}{8}}}+{q}^{{\frac{361}{8}}}+O  \left( {q}^{{\frac{401}{8}}} \right)  
$$
\emapleout
\bmaple
$>$ \enskip  etaprodtoqseries2(ep,50);
\emaple
\bmapleout
$$
(1+q+{q}^{3}+{q}^{6}+{q}^{10}+{q}^{15}+{q}^{21}+{q}^{28}+{q}^{36}+{q}^{ 45}+O \left( {q}^{50} \right) ) 
$$
\emapleout
\emaplesess
}}\vskip0pt\noindent
%%%%%%%%%%%%%%%%%%%%%%%%%%%%%%%%%%%%%%%%%%%%%%%%%%%%%%%%%
%%     |\^/|     Maple 2018 (X86 64 WINDOWS)
%% ._|\|   |/|_. Copyright (c) Maplesoft, a division of Waterloo Maple Inc. 2018
%%  \  MAPLE  /  All rights reserved. Maple is a trademark of
%%  <____ ____>  Waterloo Maple Inc.
%%       |       Type ? for help.
%% >  with(qseries):
%% >  with(ETA):
%% >  gp:=[2,2,1,-1]:
%% >  ep:=gp2etaprod(gp);
%%                                                2
%%                                      eta(2 tau)
%%                                ep := -----------
%%                                       eta(tau)
%% 
%% >  etaprodtoqseries(ep,50);
%%  (1/8)    (9/8)    (25/8)    (49/8)    (81/8)    (121/8)    (169/8)    (225/8)
%% q      + q      + q       + q       + q       + q        + q        + q
%% 
%%         (289/8)    (361/8)      (401/8)
%%      + q        + q        + O(q       )
%% 
%% >  etaprodtoqseries2(ep,50);
%%                    3    6    10    15    21    28    36    45      50
%%           1 + q + q  + q  + q   + q   + q   + q   + q   + q   + O(q  )
%% 
%% > quit
%% memory used=4.1MB, alloc=8.3MB, time=0.11
% 
% end prog2.tex
%%%%%%%%%%%%%%%%%%%%%%%%%%%%%%%%%%%%%%%%%%%%%%%%%%%%%%%%%

%%%%%%%%%%%%%%%%%%%%%%%%%%%%%%%%%%%%
\subsection{Checking modularity} 
\mylabel{subsec:modcheck}

Newman \cite{Ne59} has found necessary and sufficient conditions
under which an eta-product is a modular function on $\Gamma_0(N)$.
Let $N>0$ be a fixed integer. Here an eta-product takes the form
\begin{equation}
f(\tau) = \prod_{d\mid N} \eta(d\tau)^{m_{d}},
\mylabel{eq:etapdef}
\end{equation}
where each $d>0$ and $m_{d}\in\Z$.

\begin{theorem}[Theorem 4.7, \cite{Ne59}]
\mylabel{thm:etamodthm}
The function $f(\tau)$ (given in \eqn{etapdef}) is a modular function
on $\Gamma_0(N)$ if and only if
\begin{enumerate}
\item
$\displaystyle\sum_{d\mid N} m_d = 0$,
\item
$\displaystyle\sum_{d\mid N} d m_d \equiv0\pmod{24}$,
\item
$\displaystyle\sum_{d\mid N} \frac{N m_d}{d} \equiv0\pmod{24}$, and
\item
$\displaystyle\prod_{d\mid N} d^{|m_d|}$ is a square.
\end{enumerate}
\end{theorem}

By this theorem the eta-product \eqn{etatrdef} corresponding to the
generalized permutation 
\begin{equation}
[t_1, r_1, t_2, r_2, \dots t_m, r_m].
\mylabel{eq:gptr}
\end{equation}
is a modular function on $\Gamma_0(N)$ if the following
five conditions hold:
\begin{enumerate}
\item[(1)]
$\displaystyle \sum_{j=1}^m r_j = 0.$
\item[(2)]
$\displaystyle  \sum_{j=1}^m t_j \, r_j \equiv 0 \pmod{24}.$
\item[(3)]
The integer
$\displaystyle \prod t_j^{\abs{r_j}} $
is a square.
\item[(4)]
For each $j$, $r_j\ne0$ and $t_j \mid N$.
\item[(5)]
$\displaystyle \sum_{j=1}^m \frac{N}{t_j} \, r_j \equiv 0 \pmod{24}.$
\end{enumerate}

\noindent
\texttt{gammacheck(gp,N)} --- checks whether the eta-product corresponding 
to the generalized permutation \texttt{gp} is a modular function
on $\Gamma_0(N)$. To see a test of each condition set the following global variable
\texttt{xprint := true}.

%%%%%%%%%%%%%%%%%%%%%%%%%%%%%%%%%%%%%%%%%%%%%%%%%%%%%%%%%
\showtfile{prog3.tex}
\exampleshade{
\sn
{\tt
\bmaplesess
\bmaple
$>$ \enskip  with(qseries):  \newline
\emaple
\bmaple
$>$ \enskip  with(ETA): \newline
\emaple
\bmaple
$>$ \enskip  gp:=[1,2,2,-1,10,1,5,-2];
\emaple
\bmapleout
$$
[1,2,2,-1,10,1,5,-2] 
$$
\emapleout
\bmaple
$>$ \enskip  ep:=gp2etaprod(gp);
\emaple
\bmapleout
$$
{\frac { \left( \eta \left( \tau \right)  \right) ^{2}\eta \left( 10\, \tau \right) }{\eta \left( 2\,\tau \right)  \left( \eta \left( 5\,\tau  \right)  \right) ^{2}}} 
$$
\emapleout
\bmaple
$>$ \enskip  gammacheck(gp,10);
\emaple
\bmapleout
$$
0 
$$
\emapleout
\bmaple
$>$ \enskip  ep2:=ep\carrot2;
\emaple
\bmapleout
$$
{\frac { \left( \eta \left( \tau \right)  \right) ^{4} \left( \eta  \left( 10\,\tau \right)  \right) ^{2}}{ \left( \eta \left( 2\,\tau  \right)  \right) ^{2} \left( \eta \left( 5\,\tau \right)  \right) ^{4} }} 
$$
\emapleout
\bmaple
$>$ \enskip  gp2:=GPmake(ep2);
\emaple
\bmapleout
$$
[1,4,2,-2,10,2,5,-4] 
$$
\emapleout
\bmaple
$>$ \enskip  gammacheck(gp2,10);
\emaple
\bmapleout
$$
1 
$$
\emapleout
\emaplesess
}\vskip0pt\noindent
}
%%%%%%%%%%%%%%%%%%%%%%%%%%%%%%%%%%%%%%%%%%%%%%%%%%%%%%%%%
%%     |\^/|     Maple 2018 (X86 64 WINDOWS)
%% ._|\|   |/|_. Copyright (c) Maplesoft, a division of Waterloo Maple Inc. 2018
%%  \  MAPLE  /  All rights reserved. Maple is a trademark of
%%  <____ ____>  Waterloo Maple Inc.
%%       |       Type ? for help.
%% >  with(qseries): 
%% >  with(ETA):
%% >  gp:=[1,2,2,-1,10,1,5,-2];
%%                        gp := [1, 2, 2, -1, 10, 1, 5, -2]
%% 
%% >  ep:=gp2etaprod(gp);
%%                                         2
%%                                 eta(tau)  eta(10 tau)
%%                           ep := ----------------------
%%                                                      2
%%                                 eta(2 tau) eta(5 tau)
%% 
%% >  gammacheck(gp,10);
%%                                        0
%% 
%% >  ep2:=ep^2;
%%                                         4            2
%%                                 eta(tau)  eta(10 tau)
%%                          ep2 := -----------------------
%%                                           2           4
%%                                 eta(2 tau)  eta(5 tau)
%% 
%% >  gp2:=GPmake(ep2);
%%                        gp2 := [1, 4, 2, -2, 10, 2, 5, -4]
%% 
%% >  gammacheck(gp2,10);
%%                                        1
%% 
%% > quit
%% memory used=2.6MB, alloc=8.3MB, time=0.08
% 
% end prog3.tex
%%%%%%%%%%%%%%%%%%%%%%%%%%%%%%%%%%%%%%%%%%%%%%%%%%%%%%%%%
We considered two etaproducts:
$$
g_1(\tau)=
{\frac {  \eta \left( \tau \right) ^{2}\eta \left( 10\, \tau \right) }{\eta \left( 2\,\tau \right)   \eta \left( 5\,\tau  \right) ^{2}}},\qquad
g_2(\tau)=
{\frac {  \eta \left( \tau \right) ^{4}\eta \left( 10\, \tau \right)^2 }{\eta \left( 2\,\tau \right)^2   \eta \left( 5\,\tau  \right) ^{4}}}.
$$
We see that $g_1(\tau)$ is not a modular function on $\Gamma_0(10)$ but
its square $g_2(\tau)$ is a modular function on $\Gamma_0(10)$. To see
the reasons why $g_1(\tau)$ failed we set the global variable \texttt{xprint}
to \texttt{true}.

%%%%%%%%%%%%%%%%%%%%%%%%%%%%%%%%%%%%%%%%%%%%%%%%%%%%%%%%%
\showtfile{prog3a.tex}
\exampleshade{
\sn
{\tt
\bmaplesess
\bmaple
$>$ \enskip  xprint:=true: \newline
\emaple
\bmaple
$>$ \enskip  gammacheck(gp,10);
\emaple
\bmapleout
\vskip 0pt
\noindent
\begin{obeylines}
                               Condition (1) holds\newline
                               Condition (2) holds\newline
                           Condition (3) does not hold\newline
                               Condition (4) holds\newline
                           Condition (5) does not hold\newline
                            function is NOT invariant\newline
\end{obeylines}
$$
0 
$$
\emapleout
\emaplesess
}\vskip0pt\noindent
}
%%%%%%%%%%%%%%%%%%%%%%%%%%%%%%%%%%%%%%%%%%%%%%%%%%%%%%%%%
%%     |\^/|     Maple 2018 (X86 64 WINDOWS)
%% ._|\|   |/|_. Copyright (c) Maplesoft, a division of Waterloo Maple Inc. 2018
%%  \  MAPLE  /  All rights reserved. Maple is a trademark of
%%  <____ ____>  Waterloo Maple Inc.
%%       |       Type ? for help.
%% >  xprint:=true:
%% >  with(qseries): 
%% >  with(ETA):
%% >  gp:=[1,2,2,-1,10,1,5,-2];
%%                        gp := [1, 2, 2, -1, 10, 1, 5, -2]
%% 
%% >  ep:=gp2etaprod(gp);
%%                                         2
%%                                 eta(tau)  eta(10 tau)
%%                           ep := ----------------------
%%                                                      2
%%                                 eta(2 tau) eta(5 tau)
%% 
%% >  gammacheck(gp,10);
%%                               Condition (1) holds
%% 
%%                               Condition (2) holds
%% 
%%                           Condition (3) does not hold
%% 
%%                               Condition (4) holds
%% 
%%                           Condition (5) does not hold
%% 
%%                            function is NOT invariant
%% 
%%                                        0
%% 
%% > quit
%% memory used=2.5MB, alloc=8.3MB, time=0.09
% 
% end prog3a.tex
%%%%%%%%%%%%%%%%%%%%%%%%%%%%%%%%%%%%%%%%%%%%%%%%%%%%%%%%%
We see that $g_1(\tau)$ is not a modular function on $\Gamma_0(10)$
since it failed to satisfy Conditions (3) and (5).

To check whether an eta-product is a modular form with character use
the following function.

\noindent
\texttt{gamma0FORMCHECK(gp,N)} --- checks whether the eta-product corresponding 
to the generalized permutation \texttt{gp} is a form
on $\Gamma_0(N)$ with character. See %%\cite[Theorem 1.64]{Onbook04}.
\cite{Onbook04}.
%%%%%%%%%%%%%%%%%%%%%%%%%%%%%%%%%%%%%%%%%%%%%%%%%%%%%%%%%
\showtfile{prog4.tex}
\exampleshade{
\sn
{\tt
\bmaplesess
\bmaple
$>$ \enskip  with(qseries): \newline
\emaple
\bmaple
$>$ \enskip  with(ETA): \newline
\emaple
\bmaple
$>$ \enskip  gp2:=[1,4,2,4,4,-3,10,2,20,-1]: \newline
\emaple
\bmaple
$>$ \enskip  ep2:=gp2etaprod(gp2);
\emaple
\bmapleout
$$
{\frac { \left( \eta \left( \tau \right)  \right) ^{4} \left( \eta  \left( 2\,\tau \right)  \right) ^{4} \left( \eta \left( 10\,\tau  \right)  \right) ^{2}}{ \left( \eta \left( 4\,\tau \right)  \right) ^{ 3}\eta \left( 20\,\tau \right) }} 
$$
\emapleout
\bmaple
$>$ \enskip  gamma0FORMCHECK(gp2,40);
\emaple
\bmapleout
\begin{center}
[``N='',40,``weight='',3,``character='',[$2048000$,$- \left(  \left( 2  \right)  \right) ^{14} \left(  \left( 5 \right)  \right) ^{3}$,$d$]] 
\end{center}
\emapleout
\emaplesess
}\vskip0pt\noindent
}
%%%%%%%%%%%%%%%%%%%%%%%%%%%%%%%%%%%%%%%%%%%%%%%%%%%%%%%%%
%%     |\^/|     Maple 2018 (X86 64 WINDOWS)
%% ._|\|   |/|_. Copyright (c) Maplesoft, a division of Waterloo Maple Inc. 2018
%%  \  MAPLE  /  All rights reserved. Maple is a trademark of
%%  <____ ____>  Waterloo Maple Inc.
%%       |       Type ? for help.
%% >  with(qseries):
%% >  with(ETA):
%% >  gp2:=[1,4,2,4,4,-3,10,2,20,-1]:
%% >  ep2:=gp2etaprod(gp2);
%%                                   4           4            2
%%                           eta(tau)  eta(2 tau)  eta(10 tau)
%%                    ep2 := ----------------------------------
%%                                          3
%%                                eta(4 tau)  eta(20 tau)
%% 
%% >  gamma0FORMCHECK(gp2,40);
%%                                                             14     3
%%       ["N=", 40, "weight=", 3, "character=", [2048000, - (2)    (5) , d]]
%% 
%% > quit
%% memory used=2.5MB, alloc=8.3MB, time=0.08
% 
% end prog4.tex
%%%%%%%%%%%%%%%%%%%%%%%%%%%%%%%%%%%%%%%%%%%%%%%%%%%%%%%%%
This means that the function
$$
g(\tau) =
{\frac { \eta \left( \tau \right) ^{4}  \eta  \left( 2\,\tau \right)  ^{4}  \eta \left( 10\,\tau  \right)   ^{2}}{  \eta \left( 4\,\tau \right)   ^{ 3}\eta \left( 20\,\tau \right) }} 
$$
is a modular form on $\Gamma_0(40)$ of weight $3$ and character
$\chi(d) = \leg{-2048000}{d}= \leg{-20}{d}$.
%%                                                             14     3
%%       ["N=", 40, "weight=", 3, "character=", [2048000, - (2)    (5) , d]]
%% 

\subsection{Cusps}
\mylabel{subsec:cusps}

Chua and Lang \cite{Ch-La04}
have found a set of inequivalent 
cusps for $\Gamma_0(N)$. 
\begin{theorem}[\cite{Ch-La04}(p.354)]
\mylabel{thm:chualang}
Let N be a positive integer and foreach positive divisor $d$ of $N$ let
$e_d = (d,N/d)$. Then set
$$
\Delta = \cup_{d\mid N} S_d
$$
is a complete set of inequivalent cusps of $\Gamma_0(N)$ where
$$
S_d = \{ x_i/d\,:\,(x_i,i)=d,\quad 0\le x_i\le d-1,\quad x_i\not\equiv
x_j \pmod{e_d}\}.
$$
\end{theorem}

\noindent
\texttt{cuspmake(N)} --- returns a complete set of inequivalent
cusps of $\Gamma_0(N)$.

%%%%%%%%%%%%%%%%%%%%%%%%%%%%%%%%%%%%%%%%%%%%%%%%%%%%%%%%%
\showtfile{prog5.tex}
\exampleshade{
\sn
{\tt
\bmaplesess
\bmaple
$>$ \enskip  cuspmake(40);
\emaple
\bmapleout
$$
 \left\{ 0,\frac{1}{2},\frac{1}{4},\frac{1}{5},\frac{1}{8},\frac{1}{10},
\frac{1}{20},\frac{1}{40} \right\}  
$$
\emapleout
\emaplesess
}\vskip0pt\noindent
}
%%%%%%%%%%%%%%%%%%%%%%%%%%%%%%%%%%%%%%%%%%%%%%%%%%%%%%%%%
%%     |\^/|     Maple 2018 (X86 64 WINDOWS)
%% ._|\|   |/|_. Copyright (c) Maplesoft, a division of Waterloo Maple Inc. 2018
%%  \  MAPLE  /  All rights reserved. Maple is a trademark of
%%  <____ ____>  Waterloo Maple Inc.
%%       |       Type ? for help.
%% >  with(qseries):
%% >  with(ETA):
%% >  cuspmake(40);
%%                    {0, 1/2, 1/4, 1/5, 1/8, 1/10, 1/20, 1/40}
%% 
%% > quit
%% memory used=1.8MB, alloc=8.3MB, time=0.05
% 
% end prog5.tex
%%%%%%%%%%%%%%%%%%%%%%%%%%%%%%%%%%%%%%%%%%%%%%%%%%%%%%%%%
Biagioli \cite{Bi89} has found the fan width of the cusps of
$\Gamma_0(N)$.

\begin{lemma}[Lemma 4.2, \cite{Bi89}]
If $(r,s)=1$, then the fan width of $\Gamma_0(N)$ at $\frac{r}{s}$
is
$$
\kappa\left(\Gamma_0(N); \frac{r}{s}\right) = \frac{N}{(N,s^2)}.
$$
\end{lemma}

\texttt{fanwidth(r,N)} ---
   returns the width of the cusp $r$
   for the group $\Gamma_0(N)$.

%%%%%%%%%%%%%%%%%%%%%%%%%%%%%%%%%%%%%%%%%%%%%%%%%%%%%%%%%
\showtfile{prog6.tex}
\exampleshade{
\sn
{\tt
\bmaplesess
\bmaple
$>$ \enskip  with(ETA): \newline
\emaple
\bmaple
$>$ \enskip  fanwidth(1/8,40);
\emaple
\bmapleout
$$
5 
$$
\emapleout
\emaplesess
}\vskip0pt\noindent
}
%%%%%%%%%%%%%%%%%%%%%%%%%%%%%%%%%%%%%%%%%%%%%%%%%%%%%%%%%
%%     |\^/|     Maple 2018 (X86 64 WINDOWS)
%% ._|\|   |/|_. Copyright (c) Maplesoft, a division of Waterloo Maple Inc. 2018
%%  \  MAPLE  /  All rights reserved. Maple is a trademark of
%%  <____ ____>  Waterloo Maple Inc.
%%       |       Type ? for help.
%% >  with(ETA):
%% >  fanwidth(1/8,40);
%%                                        5
%% 
%% > quit
%% memory used=1.4MB, alloc=8.3MB, time=0.08
% 
% end prog6.tex
%%%%%%%%%%%%%%%%%%%%%%%%%%%%%%%%%%%%%%%%%%%%%%%%%%%%%%%%%

\subsection{Orders at cusps}
\mylabel{subsec:cuspords}

Ligozat \cite{Li75} has computed the order of an eta-product 
 at the cusps of $\Gamma_0(N)$.

\begin{theorem}[Theorem 4.8, \cite{Li75}]
\mylabel{thm:ordthm}
If the eta-product $f(\tau)$ (given in \eqn{etapdef})   is a modular function
on $\Gamma_0(N)$, then its order at the cusp $s=\frac{b}{c}$ 
(assuming $(b,c)=1$) is
\begin{equation}
\ord(f(\tau);s)=\sum_{d\mid N} \frac{(d,c)^2 m_d}{24d}.
\label{eq:ecord}
\end{equation}
\end{theorem}  

Following \cite[p.275]{Bi89}, \cite[p.91]{Ra} we consider the order of
a function $f$ with respect to a congruence subgroup $\Gamma$ at the
cusp $\zeta\in \mathbb{Q} \cup \{\infty\}$ and denote this by
\beq
\ORD(f,\zeta,\Gamma) = \kappa(\zeta,\Gamma)\,\ord(f;\zeta)
\mylabel{eq:ORD}.
\eeq

The following are functions for computing orders and invariant
orders of eta-products at cusps.

\noindent
\texttt{cuspord(etaprod,cusp)} --- computes the invariant order
at the given cusp of the given eta-product.

\noindent
\texttt{cuspORD(etaprod,N,cusp)} --- computes the order of the
given eta-product at the given cusp with respect to the
group $\Gamma_0(N)$.

\noindent
\texttt{cuspORDS(etaprod,CUSPS,N)} --- computes the order of the
given eta-product at each cusp in a set respect to the
group $\Gamma_0(N)$.

\exampleshade{
\sn
{\tt
\bmaplesess
\bmaple
$>$ \enskip  with(ETA): \newline
\emaple
\bmaple
$>$ \enskip  gp:=[20, -3, 10, 5, 5, -2, 4, 15, 2, -25, 1, 10]: \newline
\emaple
\bmaple
$>$ \enskip  ep:=gp2etaprod(gp);
\emaple
\bmapleout
$$
{\frac { \left( \eta \left( 10\,\tau \right)  \right) ^{5} \left( \eta  \left( 4\,\tau \right)  \right) ^{15} \left( \eta \left( \tau \right)   \right) ^{10}}{ \left( \eta \left( 20\,\tau \right)  \right) ^{3}  \left( \eta \left( 5\,\tau \right)  \right) ^{2} \left( \eta \left( 2 \,\tau \right)  \right) ^{25}}} 
$$
\emapleout
\bmaple
$>$ \enskip  gammacheck(gp,20);
\emaple
\bmapleout
$$
1 
$$
\emapleout
\bmaple
$>$ \enskip  cuspord(ep,1/4);
\emaple
\bmapleout
$$
\frac{4}{5}
$$
\emapleout
\bmaple
$>$ \enskip  fanwidth(1/4,20);
\emaple
\bmapleout
$$
5 
$$
\emapleout
\bmaple
$>$ \enskip  cuspORD(ep,20,1/4);
\emaple
\bmapleout
$$
4 
$$
\emapleout
\bmaple
$>$ \enskip  cuspORDS(ep,cuspmake(20),20);
\emaple
\bmapleout
$$
[[0,1],[\frac{1}{2},-5],[\frac{1}{4},4],[\frac{1}{5},0],[\frac{1}{10},0],
[\frac{1}{20},0]] 
$$
\emapleout
\emaplesess
}\vskip0pt\noindent
}
%%%%%%%%%%%%%%%%%%%%%%%%%%%%%%%%%%%%%%%%%%%%%%%%%%%%%%%%%
%%     |\^/|     Maple 2018 (X86 64 WINDOWS)
%% ._|\|   |/|_. Copyright (c) Maplesoft, a division of Waterloo Maple Inc. 2018
%%  \  MAPLE  /  All rights reserved. Maple is a trademark of
%%  <____ ____>  Waterloo Maple Inc.
%%       |       Type ? for help.
%% >  with(ETA):
%% >  gp:=[20, -3, 10, 5, 5, -2, 4, 15, 2, -25, 1, 10]:
%% >  ep:=gp2etaprod(gp);
%%                                    5           15         10
%%                         eta(10 tau)  eta(4 tau)   eta(tau)
%%                   ep := -------------------------------------
%%                                    3           2           25
%%                         eta(20 tau)  eta(5 tau)  eta(2 tau)
%% 
%% >  gammacheck(gp,20);
%%                                        1
%% 
%% >  cuspord(ep,1/4);
%%                                       4/5
%% 
%% >  fanwidth(1/4,20);
%%                                        5
%% 
%% >  cuspORD(ep,20,1/4);
%%                                        4
%% 
%% >  cuspORDS(ep,cuspmake(20),20);
%%          [[0, 1], [1/2, -5], [1/4, 4], [1/5, 0], [1/10, 0], [1/20, 0]]
%% 
%% > quit
%% memory used=2.5MB, alloc=8.3MB, time=0.09
% 
% end prog7a.tex
%%%%%%%%%%%%%%%%%%%%%%%%%%%%%%%%%%%%%%%%%%%%%%%%%%%%%%%%%
Let
$$
g(\tau) = {\frac {  \eta \left( 10\,\tau \right)   ^{5}  \eta  \left( 4\,\tau \right)   ^{15}  \eta \left( \tau \right)    ^{10}}{  \eta \left( 20\,\tau \right)   ^{3}   \eta \left( 5\,\tau \right)   ^{2}  \eta \left( 2 \,\tau \right)   ^{25}}} 
$$
We see that $g(\tau)$ is a modular function on $\Gamma_0(20)$,
$$
\ord(g(\tau);{1}/{4})=\frac{4}{5},
$$
and
$$
\ORD(g,{1}/{4},\Gamma_0(20)) 
= \kappa({1}/{4},\Gamma_0(20))\,\ord(g,{1}/{4})
= 5 \cdot \frac{4}{5} = 4.
$$

%% >  with(ETA):
%% >  gp:=[20, -3, 10, 5, 5, -2, 4, 15, 2, -25, 1, 10]:
%% >  ep:=gp2etaprod(gp);
%%                                    5           15         10
%%                         eta(10 tau)  eta(4 tau)   eta(tau)
%%                   ep := -------------------------------------
%%                                    3           2           25
%%                         eta(20 tau)  eta(5 tau)  eta(2 tau)
%% 
%% >  gammacheck(gp,20);
%%                                        1
%% 
%% >  cuspord(ep,20,1/4);
%%                                       4/5
%% 
%% >  fanwidth(1/4,20);
%%                                        5
%% 
%% >  cuspORD(ep,20,1/4);
%%                                        4
%% 
%% >  cuspORDS(ep,cuspmake(20),20);
%%          [[0, 1], [1/2, -5], [1/4, 4], [1/5, 0], [1/10, 0], [1/20, 0]]

%%SECTION 3 %%%%%%%%%%%%%%%%%%%%%%%%%%%%%%%%%%%%%%%%%%%%%%%%%%%%%%%%%%%%%%%%%%%%%%%
\section{Proving eta-function identities}
\mylabel{sec:proveetaid}

\subsection{Linear relations between eta-products} 
\mylabel{subsec:linrels}

          To prove a identity involving eta-products one needs to  basically
do the following.
\begin{enumerate}
\item[(i)]
Rewrite the identity in terms of eta-functions.
\item[(ii)]
Check that each term in the identity is a modular function on some
group $\Gamma_0(N)$.
\item[(iii)]
Determine the order at each cusp of $\Gamma_0(N)$ of each term
in the identity.
\item[(iv)]
Use the valence formula to determine up to which power of $q$ is needed to verify
the identity.
\item[(v)]
Finally prove the identity by carrying out the verification.
\end{enumerate}

In this section we explain how to carry out each of these steps in \maplep
Then we show how the whole process of proof can be automated.

Our method for proving eta-product identities depends
on
\begin{theorem}[The Valence Formula \cite{Ra}(p.98)]
\mylabel{thm:val}
Let $f\ne0$ be a modular form of weight $k$ with respect to a subgroup $\Gamma$ of finite index
in $\Gamma(1)=\SL_2(\mathbb{Z})$. Then
\beq
\ORD(f,\Gamma) = \frac{1}{12} \mu \, k,
\mylabel{eq:valform}
\eeq
where $\mu$ is the index of $\widehat{\Gamma}$ in $\widehat{\Gamma(1)}$,
$$
\ORD(f,\Gamma) := \sum_{\zeta\in R^{*}} \ORD(f,\zeta,\Gamma),
$$
$R^{*}$ is a fundamental region for $\Gamma$,
and
$\ORD(f,\zeta,\Gamma)$ is given in equation \eqn{ORD}.
\end{theorem}
\begin{remark}
For $\zeta\in\mathfrak{h}$,
$\ORD(f,\zeta,\Gamma)$ is defined in terms of 
 the invariant order $\ord(f,\zeta)$, which  is interpreted
in the usual sense. See \cite[p.91]{Ra} for details of this and the 
notation used.
\end{remark}

Since any modular function has weight $k=0$ and any eta-product has no zeros and no
poles on the upper-half plane we have
\begin{cor}
\mylabel{cor:valcor}
Let $f_1(\tau)$, $f_2(\tau)$, \dots, $f_n(\tau)$ be eta-products that
are modular functions on $\Gamma_0(N)$. Let $\mathcal{S}_N$ be a set of inequivalent
cusps for $\Gamma_0(N)$. Define the constant
\beq
B = \sum_{\substack{s\in\mathcal{S}_N\\s\ne \infty}}
        \mbox{min}
        (\left\{\ORD(f_j,s,\Gamma_0(N))\,:\, 1 \le j \le n\right\} \cup \{0\}),
\mylabel{eq:Bdef}
\eeq
and consider
\beq
g(\tau) := \alpha_1 f_1(\tau) + \alpha_2 f_2(\tau) + \cdots + \alpha_n f_n(\tau) + 1, 
\mylabel{eq:gdef}
\eeq
where each $\alpha_j\in\mathbb{C}$. Then
$$
g(\tau) \equiv 0
$$
if and only if
\beq
\ORD(g(\tau), \infty, \Gamma_0(N)) > -B.
\mylabel{eq:ORDBineq}
\eeq
\end{cor}

    To prove an alleged eta-product  identity, we first rewrite it in the form
\begin{equation} 
    \alpha_1 f_1(\tau) + \alpha_2 f_2(\tau) + \cdots + \alpha_n f_n(\tau) + 1 = 0,
\mylabel{eq:fid}
\end{equation}
where each $\alpha_i\in\C$ and each $f_i(\tau)$ is an eta-product of 
level $N$. We use the following algorithm:

        \vskip 10pt\noindent
{\it\footnotesize STEP 0}. \quad  Write the identity in the form \eqn{fid}. 

        \vskip 10pt\noindent
{\it\footnotesize STEP 1}. \quad  Use Theorem \thm{etamodthm} to check that
$f_j(\tau)$ is a modular function on $\Gamma_0(N)$ for each
$1 \le j \le n$.

        \vskip 10pt\noindent
{\it\footnotesize STEP 2}. \quad  Use Theorem \thm{chualang} to
find a set $\mathcal{S}_N$ of inequivalent cusps for $\Gamma_0(N)$ and the
fan width of each cusp.

        \vskip 10pt\noindent
{\it\footnotesize STEP 3}. \quad  Use Theorem \thm{ordthm} to
calculate the order of each eta-product
$f_j(\tau)$ at each cusp of $\Gamma_0(N)$.

        \vskip 10pt\noindent
{\it\footnotesize STEP 4}. \quad  Calculate
        $$
        B =
        \sum_{\substack{s\in\mathcal{S}_N\\s\ne \infty}}
        \mbox{min}
        (\left\{\ORD(f_j,s,\Gamma_0(N))\,:\, 1 \le j \le n\right\} \cup \{0\}).
        $$
        %%
        %%\left\{\ORD(f_j;s;\Gamma_1(N))\,:\, 1 \le j \le n\right\} \cup \{0\}
        %%

        \vskip 10pt\noindent
{\it\footnotesize STEP 5}. \quad  Show that
        $$
        \ORD(g(\tau),\infty,\Gamma_0(N)) > -B
        $$
        where
        $$
        g(\tau) = \alpha_1 f_1(\tau) + \alpha_2 f_2(\tau) +
        \cdots + \alpha_n f_n(\tau) + 1.
        $$
        Corollary \corol{valcor} then implies that $g(\tau)\equiv0$ and
        hence the eta-product identity  \eqn{fid}.

To calculate the constant $B$ we use 

\texttt{mintotORDS(L,n)} --- returns the constant $B$ in equation \eqn{Bdef}
where $L$ is the array of ORDS:
$$
L := [\ORD(f_1), \ORD(f_2), \dots, \ORD(f_n)],
$$
where
$$
\ORD(f) = [\ORD(f,\zeta_1,\Gamma_0(N)), 
\ORD(f,\zeta_2,\Gamma_0(N)), \dots,
\ORD(f,\zeta_m,\Gamma_0(N))] 
$$
and $\zeta_1$, $\zeta_2$, \dots, $\zeta_m$ are the 
inequivalent cusps of $\Gamma_0(N)$ apart from $\infty$.
Each $\ORD(f)$ is computed using \texttt{getaprodcuspORDS}.

\noindent
EXAMPLE:
As an example we prove Ramanujan's Entry 3.1:
\beq
PQ + \frac{9}{PQ} = \left(\frac{Q}{P}\right)^3 +  \left(\frac{P}{Q}\right)^3,  
\mylabel{eq:ramentry3p1}
\eeq
where
$$
P = \frac{\eta(\tau)^2}{\eta(3\tau)^2} \quad\mbox{and}\quad
Q = \frac{\eta(2\tau)^2}{\eta(6\tau)^2}.
$$

\vskip 10pt\noindent
{\it\footnotesize STEP 0}. \quad  We rewrite the identity \eqn{ramentry3p1}
in the form \eqn{fid}
%%%%%%%%%%%%%%%%%%%%%%%%%%%%%%%%%%%%%%%%%%%%%%%%%%%%%%%%%
\showtfile{prog8.tex}
\mapleshade{
\sn
{\tt
\bmaplesess
\bmaple
$>$ \enskip  gpP:=[1,2,3,-2]:  \newline
\emaple
\bmaple
$>$ \enskip  gpQ:=[2,2,6,-2]: \newline
\emaple
\bmaple
$>$ \enskip  P:=gp2etaprod(gpP);
\emaple
\bmapleout
$$
P:={\frac { \left( \eta \left( \tau \right)  \right) ^{2}}{ \left( \eta  \left( 3\,\tau \right)  \right) ^{2}}} 
$$
\emapleout
\bmaple
$>$ \enskip  Q:=gp2etaprod(gpQ);
\emaple
\bmapleout
$$
Q:={\frac { \left( \eta \left( 2\,\tau \right)  \right) ^{2}}{ \left( \eta  \left( 6\,\tau \right)  \right) ^{2}}} 
$$
\emapleout
\bmaple
$>$ \enskip  ETAid:=P*Q+9/P/Q - (Q/P)\carrot3 - (P/Q)\carrot3;
\emaple
\bmapleout
$$
ETAid:={\frac { \left( \eta \left( \tau \right)  \right) ^{2} \left( \eta  \left( 2\,\tau \right)  \right) ^{2}}{ \left( \eta \left( 3\,\tau  \right)  \right) ^{2} \left( \eta \left( 6\,\tau \right)  \right) ^{2} }}+9\,{\frac { \left( \eta \left( 3\,\tau \right)  \right) ^{2} \left(  \eta \left( 6\,\tau \right)  \right) ^{2}}{ \left( \eta \left( \tau  \right)  \right) ^{2} \left( \eta \left( 2\,\tau \right)  \right) ^{2} }}-{\frac { \left( \eta \left( 2\,\tau \right)  \right) ^{6} \left(  \eta \left( 3\,\tau \right)  \right) ^{6}}{ \left( \eta \left( 6\,\tau  \right)  \right) ^{6} \left( \eta \left( \tau \right)  \right) ^{6}}}- {\frac { \left( \eta \left( 6\,\tau \right)  \right) ^{6} \left( \eta  \left( \tau \right)  \right) ^{6}}{ \left( \eta \left( 2\,\tau  \right)  \right) ^{6} \left( \eta \left( 3\,\tau \right)  \right) ^{6} }} 
$$
\emapleout
\bmaple
$>$ \enskip  ETAidn:=etanormalid(\%);
\emaple
\bmapleout
$$
ETAidn := 1+9\,{\frac { \left( \eta \left( 3\,\tau \right)  \right) ^{4} \left(  \eta \left( 6\,\tau \right)  \right) ^{4}}{ \left( \eta \left( \tau  \right)  \right) ^{4} \left( \eta \left( 2\,\tau \right)  \right) ^{4} }}-{\frac { \left( \eta \left( 3\,\tau \right)  \right) ^{8} \left(  \eta \left( 2\,\tau \right)  \right) ^{4}}{ \left( \eta \left( \tau  \right)  \right) ^{8} \left( \eta \left( 6\,\tau \right)  \right) ^{4} }}-{\frac { \left( \eta \left( \tau \right)  \right) ^{4} \left( \eta  \left( 6\,\tau \right)  \right) ^{8}}{ \left( \eta \left( 3\,\tau  \right)  \right) ^{4} \left( \eta \left( 2\,\tau \right)  \right) ^{8} }} 
$$
\emapleout
\emaplesess
}\vskip0pt\noindent
}
Thus identity \eqn{ramentry3p1} is equivalent to
$$
g(\tau) = 0,
$$
where
\beq
g(\tau) = 1+9\,f_1(\tau) - f_2(\tau) - f_3(\tau),
\mylabel{eq:gdefeg}
\eeq
$$
f_1={\frac { \left( \eta \left( 3\,\tau \right)  \right) ^{4} \left( \eta  \left( 6\,\tau \right)  \right) ^{4}}{ \left( \eta \left( \tau  \right)  \right) ^{4} \left( \eta \left( 2\,\tau \right)  \right) ^{4} }} 
\quad
f_2={\frac { \left( \eta \left( 3\,\tau \right)  \right) ^{8} \left( \eta  \left( 2\,\tau \right)  \right) ^{4}}{ \left( \eta \left( \tau  \right)  \right) ^{8} \left( \eta \left( 6\,\tau \right)  \right) ^{4} }} ,
f_3={\frac { \left( \eta \left( \tau \right)  \right) ^{4} \left( \eta  \left( 6\,\tau \right)  \right) ^{8}}{ \left( \eta \left( 3\,\tau  \right)  \right) ^{4} \left( \eta \left( 2\,\tau \right)  \right) ^{8} }} .
$$

\vskip 10pt\noindent
{\it\footnotesize STEP 1}. \quad  We check that each eta-product is a modular 
function on $\Gamma_0(6)$.
%%%%%%%%%%%%%%%%%%%%%%%%%%%%%%%%%%%%%%%%%%%%%%%%%%%%%%%%%
\showtfile{prog9.tex}
\mapleshade{
\sn
{\tt
\bmaplesess
\bmaple
$>$ \enskip  f1:=op(2,ETAidn)/9;
\emaple
\bmapleout
$$
f1:={\frac { \left( \eta \left( 3\,\tau \right)  \right) ^{4} \left( \eta  \left( 6\,\tau \right)  \right) ^{4}}{ \left( \eta \left( \tau  \right)  \right) ^{4} \left( \eta \left( 2\,\tau \right)  \right) ^{4} }} 
$$
\emapleout
\bmaple
$>$ \enskip  f2:=-op(3,ETAidn);
\emaple
\bmapleout
$$
f2:={\frac { \left( \eta \left( 3\,\tau \right)  \right) ^{8} \left( \eta  \left( 2\,\tau \right)  \right) ^{4}}{ \left( \eta \left( \tau  \right)  \right) ^{8} \left( \eta \left( 6\,\tau \right)  \right) ^{4} }} 
$$
\emapleout
\bmaple
$>$ \enskip  f3:=-op(4,ETAidn);
\emaple
\bmapleout
$$
f3:={\frac { \left( \eta \left( \tau \right)  \right) ^{4} \left( \eta  \left( 6\,\tau \right)  \right) ^{8}}{ \left( \eta \left( 3\,\tau  \right)  \right) ^{4} \left( \eta \left( 2\,\tau \right)  \right) ^{8} }} 
$$
\emapleout
\bmaple
$>$ \enskip  gpf1:=GPmake(f1): gpf2:=GPmake(f2): gpf3:=GPmake(f3): \newline
\emaple
\bmaple
$>$ \enskip  gammacheck(gpf1,6),gammacheck(gpf2,6),gammacheck(gpf3,6);
\emaple
\bmapleout
$$
1,\quad 1,\quad 1                                        
$$
\emapleout
\emaplesess
}\vskip0pt\noindent
}

\vskip 10pt\noindent
{\it\footnotesize STEP 2}. \quad  We find a set of inequivalent cusps for 
$\Gamma_0(6)$ and their fan widths.
%%%%%%%%%%%%%%%%%%%%%%%%%%%%%%%%%%%%%%%%%%%%%%%%%%%%%%%%%
\showtfile{prog10.tex}
\mapleshade{
\sn
{\tt
\bmaplesess
\bmaple
$>$ \enskip  C6:=cuspmake(6);
\emaple
\bmapleout
$$
C6 :=  \left\{ 0,\frac{1}{2},1/3,\frac{1}{6} \right\}  
$$
\emapleout
\bmaple
$>$ \enskip  seq([cusp,fanwidth(cusp,6)], cusp in C6);
\emaple
\bmapleout
$$
[0,6] [\frac{1}{2},3] [\frac{1}{3},2] [\frac{1}{6},1]                                       
$$
\emapleout
\emaplesess
}\vskip0pt\noindent
}
%%%%%%%%%%%%%%%%%%%%%%%%%%%%%%%%%%%%%%%%%%%%%%%%%%%%%%%%%
%%     |\^/|     Maple 2018 (X86 64 WINDOWS)
%% ._|\|   |/|_. Copyright (c) Maplesoft, a division of Waterloo Maple Inc. 2018
%%  \  MAPLE  /  All rights reserved. Maple is a trademark of
%%  <____ ____>  Waterloo Maple Inc.
%%       |       Type ? for help.
%% >  with(ETA):
%% >  C6:=cuspmake(6);
%%                             C6 := {0, 1/2, 1/3, 1/6}
%% 
%% >  seq([cusp,fanwidth(cusp,6)], cusp in C6);
%%                       [0, 6], [1/2, 3], [1/3, 2], [1/6, 1]
%% 
%% > quit
%% memory used=1.5MB, alloc=8.3MB, time=0.08
% 
% end prog10.tex
%%%%%%%%%%%%%%%%%%%%%%%%%%%%%%%%%%%%%%%%%%%%%%%%%%%%%%%%%

\vskip 10pt\noindent
{\it\footnotesize STEP 3}. \quad  We compute $\ORD(f_j,\zeta,\Gamma_0(6))$
for each $j$ and each cusp $\zeta$ of $\Gamma_0(6)$ apart from $\infty$.
%%%%%%%%%%%%%%%%%%%%%%%%%%%%%%%%%%%%%%%%%%%%%%%%%%%%%%%%%
\showtfile{prog11.tex}
\mapleshade{
\sn
{\tt
\bmaplesess
\bmaple
$>$ \enskip  C6:=cuspmake(6) minus \{1/6\};
\emaple
\bmapleout
$$
C6 :=  \left\{ 0,{1}/{2},1/3 \right\}  
$$
\emapleout
\bmaple
$>$ \enskip  ORDS0:=cuspORDSnotoo(1,C6,6);
\emaple
\bmapleout
$$
ORDS0 := [[0,0],\quad[{1}/{2},0],\quad[{1}/{3},0]] 
$$
\emapleout
\bmaple
$>$ \enskip  ORDS1:=cuspORDSnotoo(f1,C6,6);
\emaple
\bmapleout
$$
ORDS1:= [[0,-1],\quad[{1}/{2},-1],\quad[{1}/{3},1]] 
$$
\emapleout
\bmaple
$>$ \enskip  ORDS2:=cuspORDSnotoo(f2,C6,6);
\emaple
\bmapleout
$$
ORDS2 := [[0,-1],\quad[{1}/{2},0],\quad[{1}/{3},1]] 
$$
\emapleout
\bmaple
$>$ \enskip  ORDS3:=cuspORDSnotoo(f3,C6,6);
\emaple
\bmapleout
$$
ORDS3 := [[0,0],\quad[{1}/{2},-1],\quad[{1}/{3},0]] 
$$
\emapleout
\emaplesess
}\vskip0pt\noindent
}

\vskip 10pt\noindent
{\it\footnotesize STEP 4}. \quad  We calculate the constant $B$ in \eqn{Bdef}.
%%%%%%%%%%%%%%%%%%%%%%%%%%%%%%%%%%%%%%%%%%%%%%%%%%%%%%%%%
\showtfile{prog12.tex}
\mapleshade{
\sn
{\tt
\bmaplesess
\bmaple
$>$ \enskip  mintotGAMMA0ORDS([ORDS0,ORDS1,ORDS2,ORDS3],4);
\emaple
\bmapleout
$$
-2 
$$
\emapleout
\emaplesess
}\vskip0pt\noindent
}

\vskip 10pt\noindent
{\it\footnotesize STEP 5}. \quad  To prove the identity \eqn{ramentry3p1} 
we need to
verify that 
$$
\ORD(g(\tau),\infty,\Gamma_0(6)) > 2.
$$
%%%%%%%%%%%%%%%%%%%%%%%%%%%%%%%%%%%%%%%%%%%%%%%%%%%%%%%%%
\showtfile{prog13.tex}
\mapleshade{
\sn
{\tt
\bmaplesess
\bmaple
$>$ \enskip  qetacombo(1+9*f1-f2-f3,100);
\emaple
\bmapleout
$$
0 
$$
\emapleout
\emaplesess
}\vskip0pt\noindent
}
%%%%%%%%%%%%%%%%%%%%%%%%%%%%%%%%%%%%%%%%%%%%%%%%%%%%%%%%%
%%     |\^/|     Maple 2018 (X86 64 WINDOWS)
%% ._|\|   |/|_. Copyright (c) Maplesoft, a division of Waterloo Maple Inc. 2018
%%  \  MAPLE  /  All rights reserved. Maple is a trademark of
%%  <____ ____>  Waterloo Maple Inc.
%%       |       Type ? for help.
%% >  with(ETA):
%% >  gpP:=[1,2,3,-2]: gpQ:=[2,2,6,-2]:
%% >  P:=gp2etaprod(gpP):
%% >  Q:=gp2etaprod(gpQ):
%% >  ETAid:=P*Q+9/P/Q - (Q/P)^3 - (P/Q)^3:
%% >  ETAidn:=etanormalid(%):
%% >  f1:=op(2,ETAidn)/9:
%% >  f2:=-op(3,ETAidn):
%% >  f3:=-op(4,ETAidn):
%% >  qetacombo(1+9*f1-f2-f3,100);
%%                                        0
%% 
%% > quit
%% memory used=4.4MB, alloc=8.3MB, time=0.06
% 
% end prog13.tex
%%%%%%%%%%%%%%%%%%%%%%%%%%%%%%%%%%%%%%%%%%%%%%%%%%%%%%%%%

This completes the proof of the identity \eqn{ramentry3p1}.
We only had to show that the coefficient of 
$q^j$ was zero in the $q$-expansion of 
$g(\tau)$ for $j \le 3$.  We actually did it for 
$j \le 100$ as a check.

STEPS 1--5 may be automated using the following function.

\medskip

\texttt{provemodfuncGAMMA0id(etaid,N)} --- 
returns the constant $B$ in equation \eqn{Bdef}
and prints details of the verification and proof of the identity corresponding
to \texttt{etaid}, 
which is a linear combination of symbolic eta-products, 
and $N$ is the level. If \texttt{xprint=true} then more details of the
verification are printed. When this function is called there is a query asking
whether to verify the identity. Enter \texttt{yes} to carry out the verification.
%%%%%%%%%%%%%%%%%%%%%%%%%%%%%%%%%%%%%%%%%%%%%%%%%%%%%%%%%
\showtfile{prog13.tex}
\mapleshade{
\sn
{\tt
\bmaplesess
\bmaple
$>$ \enskip  provemodfuncGAMMA0id(1+9*f1-f2-f3,6);\newline
\emaple
\bmapleout
"TERM ", 1, "of ", 4, " *****************"\newline
"TERM ", 2, "of ", 4, " *****************"\newline
"TERM ", 3, "of ", 4, " *****************"\newline
"TERM ", 4, "of ", 4, " *****************"\newline
"mintotord = ", -2\newline
"TO PROVE the identity we need to show that v[oo](ID) > ", 2\newline
*** There were NO errors. \newline
*** o Each term was modular function on\newline
      Gamma0(6). \newline
*** o We also checked that the total order of\newline
      each term was zero.\newline
"*** WARNING: some terms were constants. ***"\newline
"See array CONTERMS."\newline
To prove the identity we will need to verify if up to \newline
q\carrot(3).\newline
Do you want to prove the identity? (yes/no)\newline
You entered yes.\newline
We verify the identity to O(q\carrot(14)).\newline
RESULT: The identity holds to O(q\carrot(14)).\newline
CONCLUSION: This proves the identity since we had only\newline
            to show that v[oo](ID) > 2.\newline
\emapleout
\emaplesess
}\vskip0pt\noindent
}

%%\texttt{provemodfuncidBATCH(JACID,N)} --- is a version of 
%%\texttt{provemodfuncid} that prints less detail and does not query.

%%\exampleshade{
%%> provemodfuncidBATCH(JACID,25);\\
%%*** There were NO errors.  Each term was modular function on\\
%%    Gamma1(25). Also -mintotord=9. To prove the identity\\
%%    we need to  check up to O(q\carrot(11)).\\
%%    To be on the safe side we check up to O(q\carrot(59)).\\
%%*** The identity is PROVED!\\
%%}

\texttt{printETAIDORDStable} --- prints an ORDs table for the $f_j$
and lower bound for $g$ after \texttt{provemodfuncGAMMA0id} is run.
Formatted output from our example is given below in Table \ref{tab:ordfs}. 
By summing the last column
we see that $B=-2$, which confirms an earlier calculation using 
\texttt{mintotORDS}. 

$$
\begin{array}{c|c|c|c|c}
\noalign{\hrule}
 \zeta&\ORD \left(f_1,\zeta\right) &\ORD \left(f_2,\zeta\right) &
\ORD \left(f_3,\zeta\right) &\mbox{Lower bound for $\ORD \left(g,\zeta\right)$} 
\\ 
\noalign{\hrule}
%%{\it oo}&-1&-1&1&-1\\ 0&1&0&0&0 \\ 
0& -1&-1&0&-1\\ 
\frac{1}{2}&-1&0&-1&-1\\ 
\frac{1}{3}&1&1&0&0\\
\noalign{\hrule}
\end{array}
$$
{Orders at the cusps of $\Gamma_0(6)$ of the functions 
$f_1$, $f_2$, $f_3$ and $g$  in 
\eqn{gdefeg} needed in the proof of Ramanujan's identity \eqn{ramentry3p1}.
This table was produced by
\texttt{printETAIDORDStable()}.}                                              
\mylabel{tab:ordfs}

\texttt{provemodfuncGAMMA0idBATCH(etaid,N)} --- is a version of 
\texttt{provemodfuncGAMMA0id} that prints less detail and does not query.

%%%%%%%%%%%%%%%%%%%%%%%%%%%%%%%%%%%%%%%%%%%%%%%%%%%%%%%%%
\showtfile{prog15.tex}
\mapleshade{
\sn
{\tt
\bmaplesess
\bmaple
$>$ \enskip  provemodfuncGAMMA0idBATCH(1+9*f1-f2-f3,6);\newline
\emaple
\bmapleout
 *** There were NO errors. \newline
 *** o Each term was modular function on\newline
       Gamma0(6). \newline
 *** o We also checked that the total order of\newline
       each term was zero.\newline
 To prove the identity we will need to verify if up to \newline
 q\carrot(3).\newline
 *** The identity below is PROVED!\newline
%%                                     [1, -2, 0]
%% 
$$
[1,-2,0] 
$$
\emapleout
\emaplesess
}\vskip0pt\noindent
}
%%%%%%%%%%%%%%%%%%%%%%%%%%%%%%%%%%%%%%%%%%%%%%%%%%%%%%%%%
%%     |\^/|     Maple 2018 (X86 64 WINDOWS)
%% ._|\|   |/|_. Copyright (c) Maplesoft, a division of Waterloo Maple Inc. 2018
%%  \  MAPLE  /  All rights reserved. Maple is a trademark of
%%  <____ ____>  Waterloo Maple Inc.
%%       |       Type ? for help.
%% >  with(ETA):
%% >  gpP:=[1,2,3,-2]: gpQ:=[2,2,6,-2]:
%% >  P:=gp2etaprod(gpP):
%% >  Q:=gp2etaprod(gpQ):
%% >  ETAid:=P*Q+9/P/Q - (Q/P)^3 - (P/Q)^3:
%% >  ETAidn:=etanormalid(%):
%% >  f1:=op(2,ETAidn)/9:
%% >  f2:=-op(3,ETAidn):
%% >  f3:=-op(4,ETAidn):
%% >  provemodfuncGAMMA0idBATCH(1+9*f1-f2-f3,6);
%% *** There were NO errors. 
%% *** o Each term was modular function on
%%       Gamma0(6). 
%% *** o We also checked that the total order of
%%       each term was zero.
%% To prove the identity we will need to verify if up to 
%% q^(3).
%% *** The identity below is PROVED!
%%                                     [1, -2, 0]
%% 
%% > quit
%% memory used=3.2MB, alloc=8.3MB, time=0.09
% 
% end prog15.tex
%%%%%%%%%%%%%%%%%%%%%%%%%%%%%%%%%%%%%%%%%%%%%%%%%%%%%%%%%

Let $L=[1,-2]$. $L[1]=1$ means the identity is proved. $L[2]=B$ (the constant
in equation \eqn{Bdef}), and  we see that $B=-2$. 
This confirms an earlier calculation using 
\texttt{mintotORDS}. To print out minimal information set \texttt{noprint:=true}.

%%%%%%%%%%%%%%%%%%%%%%%%%%%%%%%%%%%%%%%%%%%%%%%%%%%%%%%%%
\showtfile{prog16.tex}
\mapleshade{
\sn
{\tt
\bmaplesess
\bmaple
$>$ \enskip  noprint:=true: \newline
\emaple
\bmaple
$>$ \enskip  provemodfuncGAMMA0idBATCH(1+9*f1-f2-f3,6);
\emaple
\bmapleout
$$
[1,-2,0] 
$$
\emapleout
\emaplesess
}\vskip0pt\noindent
}
%%%%%%%%%%%%%%%%%%%%%%%%%%%%%%%%%%%%%%%%%%%%%%%%%%%%%%%%%
%%     |\^/|     Maple 2018 (X86 64 WINDOWS)
%% ._|\|   |/|_. Copyright (c) Maplesoft, a division of Waterloo Maple Inc. 2018
%%  \  MAPLE  /  All rights reserved. Maple is a trademark of
%%  <____ ____>  Waterloo Maple Inc.
%%       |       Type ? for help.
%% >  with(ETA):
%% >  gpP:=[1,2,3,-2]: gpQ:=[2,2,6,-2]:
%% >  P:=gp2etaprod(gpP):
%% >  Q:=gp2etaprod(gpQ):
%% >  ETAid:=P*Q+9/P/Q - (Q/P)^3 - (P/Q)^3:
%% >  ETAidn:=etanormalid(%):
%% >  f1:=op(2,ETAidn)/9:
%% >  f2:=-op(3,ETAidn):
%% >  f3:=-op(4,ETAidn):
%% >  noprint:=true:
%% >  provemodfuncGAMMA0idBATCH(1+9*f1-f2-f3,6);
%%                                     [1, -2, 0]
%% 
%% > quit
%% memory used=3.1MB, alloc=8.3MB, time=0.09
% 
% end prog16.tex
%%%%%%%%%%%%%%%%%%%%%%%%%%%%%%%%%%%%%%%%%%%%%%%%%%%%%%%%%

\subsection{$\Up_p$ identities} 
\mylabel{subsec:Upids}

Let $p$ be prime and suppose
$$
f(\tau) = \sum_{n\ge n_0} a(n) q^n 
$$
is a modular function.
We define 
\begin{equation}
\Up_p f(\tau) = \sum_{pn\ge n_0} a(pn) q^n = 
 \frac{1}{{p}} \sum_{k=0}^{p-1} f\Lpar{\frac{\tau+k}{p}}.
\mylabel{eq:Updef}
\end{equation}
In this section show how to modify out method to prove eta-product
identities for $\Up_{p} f(\tau)$ when $f(\tau)$ is also an eta-product.

It is known that if $f(\tau)$ is a modular function on $\Gamma_0(pN)$,
where $p\mid N$, then $\Up_{p} f(\tau)$ is a modular function on
$\Gamma_0(N)$. 
Gordon and Hughes \cite[Theorem 4, p.336]{GoHu81} have found
lower bounds for the invariant orders of $\Up_{p} f(\tau)$ at cusps.
Let $\nu_p(n)$ denote the $p$-adic order of an integer $n$; i.e.\ the
highest power of $p$ that divides $n$.

\begin{theorem}[Theorem 4, \cite{GoHu81}]
\mylabel{thm:UpLB}
Suppose $f(\tau)$ is a modular function on $\Gamma_0(pN)$, where $p$ is
prime and $p\mid N$. Let $r=\frac{\beta}{\delta}$ be a cusp of
$\Gamma_0(N)$, where $\delta\mid N$ and $(\beta,\delta)=1$. Then
$$
\ORD(\Up_{p} f, r, \Gamma_0(N)) \ge 
\begin{cases}
\frac{1}{p} \ORD(f,r/p,\Gamma_0(pN)) & \mbox{if $\nu_p(\delta)\ge \frac{1}{2} \nu_p(N)$}\\
\ORD(f,r/p,\Gamma_0(pN)) & \mbox{if $0 < \nu_p(\delta)< \frac{1}{2} \nu_p(N)$}\\
\umin{0 \le k \le p-1} \ORD(f,(r+k)/p,\Gamma_0(pN)) 
&\mbox{if $\nu_p(\delta)=0$}.
\end{cases}
$$
\end{theorem}
%%%%%%%%%%%%%%%%%%%%%%%%%%%%%%%%%%%%%%%%%%%%%%%%%%%%%%%%%
\showtfile{prog17.tex}
\exampleshade{
\sn
{\tt
\bmaplesess
\bmaple
$>$ \enskip  with(qseries): \newline
\emaple
\bmaple
$>$ \enskip  with(ETA): \newline
\emaple
\bmaple
$>$ \enskip  gpF:=[2,1,25,1,1,-1,50,-1]: \newline
\emaple
\bmaple
$>$ \enskip  epF:=gp2etaprod(gpF);
\emaple
\bmapleout
$$
epF:={\frac {\eta \left( 2\,\tau \right) \eta \left( 25\,\tau \right) }{\eta  \left( \tau \right) \eta \left( 50\,\tau \right) }} 
$$
\emapleout
\bmaple
$>$ \enskip  gammacheck(gpF,50);
\emaple
\bmapleout
$$
1 
$$
\emapleout
\bmaple
$>$ \enskip  seq([cusp,UpLB(epF,cusp,50,5)],cusp in cuspmake(10));
\emaple
\bmapleout
$$
[0,-1],\quad  [{1}/{2},0],\quad  [{1}/{5},1/5],\quad  [1/10,-{1}/{5}]
$$
\emapleout
\emaplesess
}\vskip0pt\noindent
}
%%%%%%%%%%%%%%%%%%%%%%%%%%%%%%%%%%%%%%%%%%%%%%%%%%%%%%%%%
%%     |\^/|     Maple 2018 (X86 64 WINDOWS)
%% ._|\|   |/|_. Copyright (c) Maplesoft, a division of Waterloo Maple Inc. 2018
%%  \  MAPLE  /  All rights reserved. Maple is a trademark of
%%  <____ ____>  Waterloo Maple Inc.
%%       |       Type ? for help.
%% >  with(qseries):
%% >  with(ETA):
%% >  gpF:=[2,1,25,1,1,-1,50,-1]:
%% >  epF:=gp2etaprod(gpF);
%%                                 eta(2 tau) eta(25 tau)
%%                          epF := ----------------------
%%                                  eta(tau) eta(50 tau)
%% 
%% >  gammacheck(gpF,50);
%%                                        1
%% 
%% >  epF:=gp2etaprod(gpF);
%%                                 eta(2 tau) eta(25 tau)
%%                          epF := ----------------------
%%                                  eta(tau) eta(50 tau)
%% 
%% >  seq([cusp,UpLB(epF,cusp,50,5)],cusp in cuspmake(10));
%%                   [0, -1], [1/2, 0], [1/5, 1/5], [1/10, -1/5]
%% 
%% > quit
%% memory used=2.7MB, alloc=8.3MB, time=0.11
% 
% end prog17.tex
%%%%%%%%%%%%%%%%%%%%%%%%%%%%%%%%%%%%%%%%%%%%%%%%%%%%%%%%%
We see that
$$
F(\tau)={\frac {\eta \left( 2\,\tau \right) \eta \left( 25\,\tau \right) }{\eta  \left( \tau \right) \eta \left( 50\,\tau \right) }} 
$$
is a modular function on $\Gamma_0(50)$. From Theorem \thm{UpLB} we
have the following lower bounds for the orders of $\Up_5 F$
at the cusps of $\Gamma_0(10)$.
%%                   [0, -1], [1/2, 0], [1/5, 1/5], [1/10, -1/5]
$$
\begin{array}{|c|c|c|c|c|}
\noalign{\hrule}
\zeta   &  0  &  1/2  &   1/5  &  1/10 \\
\noalign{\hrule}
\ORD\left(\Up_5 F,\zeta, \Gamma_0(10)\right) \ge &
          -1  &  0    &   1/5  &  -1/5\\
\noalign{\hrule}
\end{array}
$$
This example is taken from \cite[p.338]{GoHu81}. It turns out
that $\Up_5 F$ is an etaproduct.
%%%%%%%%%%%%%%%%%%%%%%%%%%%%%%%%%%%%%%%%%%%%%%%%%%%%%%%%%
\showtfile{prog18.tex}
\mapleshade{
\sn
{\tt
\bmaple
$>$ \enskip  F:=etaprodtoqseries(epF,1000): \newline
\emaple
\bmaple
$>$ \enskip  sf:=sift(F,q,5,0,1000): \newline
\emaple
\bmaple
$>$ \enskip  epG:=etamake(sf,q,100);
\emaple
\bmapleout
$$
epG := {\frac { \left( \eta \left( 5\,\tau \right)  \right) ^{4} \left( \eta  \left( 2\,\tau \right)  \right) ^{2}}{ \left( \eta \left( 10\,\tau  \right)  \right) ^{2} \left( \eta \left( \tau \right)  \right) ^{4}}} 
$$
\emapleout
\bmaple
$>$ \enskip  seq([cusp,cuspORD(epG,10,cusp)],cusp in cuspmake(10));
\emaple
\bmapleout
$$
[0,-1],\quad  [{1}/{2},0],\quad  [{1}/{5},1],\quad  [1/10,0]
$$
\emapleout
\emaplesess
}\vskip0pt\noindent
}
%%%%%%%%%%%%%%%%%%%%%%%%%%%%%%%%%%%%%%%%%%%%%%%%%%%%%%%%%
%%     |\^/|     Maple 2018 (X86 64 WINDOWS)
%% ._|\|   |/|_. Copyright (c) Maplesoft, a division of Waterloo Maple Inc. 2018
%%  \  MAPLE  /  All rights reserved. Maple is a trademark of
%%  <____ ____>  Waterloo Maple Inc.
%%       |       Type ? for help.
%% >  with(qseries):
%% >  with(ETA):
%% >  gpF:=[2,1,25,1,1,-1,50,-1]:
%% >  epF:=gp2etaprod(gpF);
%%                                 eta(2 tau) eta(25 tau)
%%                          epF := ----------------------
%%                                  eta(tau) eta(50 tau)
%% 
%% >  F:=etaprodtoqseries(epF,1000):
%% memory used=15.8MB, alloc=43.6MB, time=0.14
%% >  sf:=sift(F,q,5,0,1000):
%% >  epG:=etamake(sf,q,100);
%%                                           4           2
%%                                 eta(5 tau)  eta(2 tau)
%%                          epG := -----------------------
%%                                            2         4
%%                                 eta(10 tau)  eta(tau)
%% 
%% >  seq([cusp,cuspORD(epG,10,cusp)],cusp in cuspmake(10));
%%                      [0, -1], [1/2, 0], [1/5, 1], [1/10, 0]
%% 
%% > quit
%% memory used=24.0MB, alloc=43.6MB, time=0.17
% 
% end prog18.tex
%%%%%%%%%%%%%%%%%%%%%%%%%%%%%%%%%%%%%%%%%%%%%%%%%%%%%%%%%
We see that
$$
\Up_5(F) = G,
$$
where
$$
G(\tau) = {\frac {  \eta \left( 5\,\tau \right)   ^{4}  \eta  \left( 2\,\tau  \right) ^{2}}{  \eta \left( 10\,\tau  \right)  ^{2}  \eta \left( \tau \right)  ^{4}}}.
$$
%%                      [0, -1], [1/2, 0], [1/5, 1], [1/10, 0]
We have the following exact values for the orders of $G(\tau)$
at the cusps of $\Gamma_0(10)$.
$$
\begin{array}{|c|c|c|c|c|}
\noalign{\hrule}
\zeta   &  0  &  1/2  &   1/5  &  1/10 \\
\noalign{\hrule}
\ORD\left(G,\zeta, \Gamma_0(10)\right) = &
          -1  &  0    &   1  &  0\\
\noalign{\hrule}
\end{array}
$$
These values are consistent with the lower bounds that we found.

It is a simple matter to modify our method, from Section \subsect{linrels},  
for proving 
linear relations between eta-products, to proving $\Up_p$ eta-product
identities. We wish to prove an identity of the form
\begin{equation} 
    \Up_p(g) = 
\alpha_1 f_1(\tau) + \alpha_2 f_2(\tau) + \cdots + \alpha_n f_n(\tau),
\mylabel{eq:Upfid}
\end{equation}
where $p$ is prime, $p\mid N$, $g(\tau)$ is an eta-product and 
a modular function on $\Gamma_0(pN)$, and each $f_j(\tau)$
is an eta-product and modular function on $\Gamma_0(N)$.
 We use the following algorithm:

        \vskip 10pt\noindent
{\it\footnotesize STEP 0}. \quad  Write the identity in the form \eqn{Upfid}. 

        \vskip 10pt\noindent
{\it\footnotesize STEP 1}. \quad  Use Theorem \thm{etamodthm} to check that
$f_j(\tau)$ is a modular function on $\Gamma_0(N)$ for each
$1 \le j \le n$, and $g(\tau)$ is a modular function on $\Gamma_0(pN)$.

        \vskip 10pt\noindent
{\it\footnotesize STEP 2}. \quad  Use Theorem \thm{chualang} to
find a set $\mathcal{S}_N$ of inequivalent cusps for $\Gamma_0(N)$ and the
fan width of each cusp.

        \vskip 10pt\noindent
{\it\footnotesize STEP 3a}. \quad  We compute $\ORD(f_j,\zeta,\Gamma_0(N))$
for each $j$ and each cusp $\zeta$ of $\Gamma_0(N)$ apart from $\infty$.

        \vskip 10pt\noindent
{\it\footnotesize STEP 3b}. \quad  Use Theorem \thm{UpLB}
to find a lower bound for
$$
\ORD(\Up_{p} g, r, \Gamma_0(N)) 
$$
for each cusp $r$ of $\Gamma_0(N)$.  Call this 
lower bound $L(g,r,N)$.

        \vskip 10pt\noindent
{\it\footnotesize STEP 4}. \quad  Calculate
\beq
        B =
        \sum_{\substack{s\in\mathcal{S}_N\\s\ne \infty}}
        \mbox{min}
(\left\{\ORD(f_j,s,\Gamma_0(N))\,:\, 1 \le j \le n\right\} \cup \{L(g,s,N)\}).
\mylabel{eq:UpB}
\eeq
        %%
        %%\left\{\ORD(f_j;s;\Gamma_1(N))\,:\, 1 \le j \le n\right\} \cup \{0\}
        %%

        \vskip 10pt\noindent
{\it\footnotesize STEP 5}. \quad  Show that
        $$
        \ORD(h(\tau),\infty,\Gamma_0(N)) > -B
        $$
        where
        $$
        h(\tau) = \Up_p(g)  
        - (\alpha_1 f_1(\tau) + \alpha_2 f_2(\tau) +
        \cdots + \alpha_n f_n(\tau)).
        $$
        Corollary \corol{valcor} then implies that $h(\tau)\equiv0$ and
        hence the $\Up_p$ eta-product identity  \eqn{Upfid}.

\noindent
EXAMPLE:
Let
$$
g(\tau)={\frac {  \eta \left( 50\,\tau \right)   ^{5}  \eta  \left( 5\,\tau\right)    ^{4}  \eta \left( 4\,\tau  \right)   ^{3}  \eta \left( 2\,\tau \right)   ^{3} }{  \eta \left( 100\,\tau \right)   ^{3}  \eta  \left( 25\,\tau \right)   ^{2}  \eta \left( 10\,\tau  \right)   ^{8}  \eta \left( \tau \right)   ^{2}}}.
$$
We prove that
\beq
\Up_5(g(\tau)) =
5\,{\frac {  \eta \left( 10\,\tau \right)   ^{8}   \eta \left( \tau \right)   ^{4}}{  \eta \left( 5\,\tau  \right)   ^{4}  \eta \left( 2\,\tau \right)   ^{8} }}+2\,{\frac {  \eta \left( 10\,\tau \right)   ^{5}   \eta \left( \tau \right)   ^{2}}{  \eta \left( 20 \,\tau \right)   ^{3}  \eta \left( 5\,\tau \right)    ^{2}\eta \left( 4\,\tau \right) \eta \left( 2\,\tau \right) }} 
\mylabel{eq:Up5id}
\eeq

        \vskip 10pt\noindent
{\it\footnotesize STEP 0}. \quad  Write the identity in the form \eqn{Upfid}. 
The identity is
$$
\Up_5(g(\tau)) = 5\,f_1(\tau) + 2\,f_2(\tau),
$$
where
$$
f_1(\tau) = {\frac {  \eta \left( 10\,\tau \right)  ^{8}  \eta  \left( \tau \right)  ^{4}}{  \eta \left( 5\,\tau  \right)  ^{4}  \eta \left( 2\,\tau \right)  ^{8} }},\quad
f_2(\tau)= {\frac {  \eta \left( 10\,\tau \right)  ^{5}  \eta  \left( \tau \right)  ^{2}}{  \eta \left( 20\,\tau  \right)  ^{3}  \eta \left( 5\,\tau \right)  ^{2} \eta \left( 4\,\tau \right) \eta \left( 2\,\tau \right) }}.
$$

        \vskip 10pt\noindent
{\it\footnotesize STEP 1}. \quad  Use Theorem \thm{etamodthm} to check that
$f_j(\tau)$ is a modular function on $\Gamma_0(20)$ for each
$1 \le j \le 2$, and $g(\tau)$ is a modular function on $\Gamma_0(100)$.
%%%%%%%%%%%%%%%%%%%%%%%%%%%%%%%%%%%%%%%%%%%%%%%%%%%%%%%%%
\showtfile{prog19.tex}
\mapleshade{
\sn
{\tt
\bmaplesess
\bmaple
$>$ \enskip  with(qseries): \newline
\emaple
\bmaple
$>$ \enskip  with(ETA): \newline
\emaple
\bmaple
$>$ \enskip  gpg:=[100, -3, 50, 5, 25, -2, 10, -8, 5, 4, 4, 3, 2, 3, 1, -2]: \newline
\emaple
\bmaple
$>$ \enskip  epg:=gp2etaprod(gpg);
\emaple
\bmapleout
$$
epg:={\frac { \left( \eta \left( 50\,\tau \right)  \right) ^{5} \left( \eta  \left( 5\,\tau \right)  \right) ^{4} \left( \eta \left( 4\,\tau  \right)  \right) ^{3} \left( \eta \left( 2\,\tau \right)  \right) ^{3} }{ \left( \eta \left( 100\,\tau \right)  \right) ^{3} \left( \eta  \left( 25\,\tau \right)  \right) ^{2} \left( \eta \left( 10\,\tau  \right)  \right) ^{8} \left( \eta \left( \tau \right)  \right) ^{2}}} 
$$
\emapleout
\bmaple
$>$ \enskip  gammacheck(gpg,100);
\emaple
\bmapleout
$$
1 
$$
\emapleout
\bmaple
$>$ \enskip  gpf1:=[10, 8, 5, -4, 2, -8, 1, 4]: \newline
\emaple
\bmaple
$>$ \enskip  epf1:=gp2etaprod(gpf1);
\emaple
\bmapleout
$$
epf1:={\frac { \left( \eta \left( 10\,\tau \right)  \right) ^{8} \left( \eta  \left( \tau \right)  \right) ^{4}}{ \left( \eta \left( 5\,\tau  \right)  \right) ^{4} \left( \eta \left( 2\,\tau \right)  \right) ^{8} }} 
$$
\emapleout
\bmaple
$>$ \enskip  gpf2:=[20, -3, 10, 5, 5, -2, 4, -1, 2, -1, 1, 2]: \newline
\emaple
\bmaple
$>$ \enskip  epf2:=gp2etaprod(gpf2);
\emaple
\bmapleout
$$
epf2:={\frac { \left( \eta \left( 10\,\tau \right)  \right) ^{5} \left( \eta  \left( \tau \right)  \right) ^{2}}{ \left( \eta \left( 20\,\tau  \right)  \right) ^{3} \left( \eta \left( 5\,\tau \right)  \right) ^{2} \eta \left( 4\,\tau \right) \eta \left( 2\,\tau \right) }} 
$$
\emapleout
\bmaple
$>$ \enskip  gammacheck(gpf1,20),gammacheck(gpf2,20);
\emaple
\bmapleout
$$
1,\quad  1
$$
\emapleout
\emaplesess
}\vskip0pt\noindent
}
        \vskip 10pt\noindent
{\it\footnotesize STEP 2}. \quad  Use Theorem \thm{chualang} to
find a set $\mathcal{S}_{20}$ of inequivalent cusps for $\Gamma_0(20)$ and the
fan width of each cusp.
%%%%%%%%%%%%%%%%%%%%%%%%%%%%%%%%%%%%%%%%%%%%%%%%%%%%%%%%%
\showtfile{prog20.tex}
\mapleshade{
\sn
{\tt
\bmaplesess
\bmaple
$>$ \enskip  with(ETA): \newline
\emaple
\bmaple
$>$ \enskip  C20:=cuspmake(20);
\emaple
\bmapleout
$$
 \left\{ 0,{1}/{2},1/4,{1}/{5},1/10,1/20 \right\}  
$$
\emapleout
\bmaple
$>$ \enskip  seq([cusp,fanwidth(cusp,20)], cusp in C20);
\emaple
\bmapleout
$$
[0,20],\quad [{1}/{2},5],\quad [{1}/{4},5],\quad [{1}/{5},4],\quad [1/10,1],\quad [1/20,1]
$$
\emapleout
\emaplesess
}\vskip0pt\noindent
}
%%%%%%%%%%%%%%%%%%%%%%%%%%%%%%%%%%%%%%%%%%%%%%%%%%%%%%%%%
%%     |\^/|     Maple 2018 (X86 64 WINDOWS)
%% ._|\|   |/|_. Copyright (c) Maplesoft, a division of Waterloo Maple Inc. 2018
%%  \  MAPLE  /  All rights reserved. Maple is a trademark of
%%  <____ ____>  Waterloo Maple Inc.
%%       |       Type ? for help.
%% >  with(ETA):
%% >  C20:=cuspmake(20);
%%                      C20 := {0, 1/2, 1/4, 1/5, 1/10, 1/20}
%% 
%% >  seq([cusp,fanwidth(cusp,20)],\quad, cusp in C20);
%%           [0, 20],\quad, [1/2, 5],\quad, [1/4, 5],\quad, [1/5, 4],\quad, [1/10, 1],\quad, [1/20, 1],\quad
%% 
%% > quit
%% memory used=1.6MB, alloc=8.3MB, time=0.05
% 
% end prog20.tex
%%%%%%%%%%%%%%%%%%%%%%%%%%%%%%%%%%%%%%%%%%%%%%%%%%%%%%%%%

        \vskip 10pt\noindent
{\it\footnotesize STEP 3a}. \quad  We compute $\ORD(f_j,\zeta,\Gamma_0(N))$
for each $j$ and each cusp $\zeta$ of $\Gamma_0(20)$ apart from $\infty$.
%%%%%%%%%%%%%%%%%%%%%%%%%%%%%%%%%%%%%%%%%%%%%%%%%%%%%%%%%
\showtfile{prog21.tex}
\mapleshade{
\sn
{\tt
\bmaple
$>$ \enskip  C20:=cuspmake(20) minus \{1/20\}: \newline
\emaple
\bmaple
$>$ \enskip  ORDS1:=cuspORDSnotoo(epf1,C20,20);
\emaple
\bmapleout
$$
ORDS1:=[[0,0],\quad[{1}/{2},-2],\quad[{1}/{4},-2],\quad[{1}/{5},0],\quad[1/10,2] 
$$
\emapleout
\bmaple
$>$ \enskip  ORDS2:=cuspORDSnotoo(epf2,C20,20);
\emaple
\bmapleout
$$
ORDS2:=[[0,1],\quad[{1}/{2},0],\quad[{1}/{4},-1],\quad[{1}/{5},0],\quad[1/10,1]]
$$
\emapleout
\emaplesess
}\vskip0pt\noindent
}
%%%%%%%%%%%%%%%%%%%%%%%%%%%%%%%%%%%%%%%%%%%%%%%%%%%%%%%%%
%%     |\^/|     Maple 2018 (X86 64 WINDOWS)
%% ._|\|   |/|_. Copyright (c) Maplesoft, a division of Waterloo Maple Inc. 2018
%%  \  MAPLE  /  All rights reserved. Maple is a trademark of
%%  <____ ____>  Waterloo Maple Inc.
%%       |       Type ? for help.
%% >  with(qseries):
%% >  with(ETA):
%% >  gpg:=[100, -3, 50, 5, 25, -2, 10, -8, 5, 4, 4, 3, 2, 3, 1, -2],\quad:
%% >  epg:=gp2etaprod(gpg):
%% >  gammacheck(gpg,100):
%% >  gpf1:=[10, 8, 5, -4, 2, -8, 1, 4],\quad:
%% >  epf1:=gp2etaprod(gpf1):
%% >  gpf2:=[20, -3, 10, 5, 5, -2, 4, -1, 2, -1, 1, 2],\quad:
%% >  epf2:=gp2etaprod(gpf2):
%% >  C20:=cuspmake(20) minus {1/20}:
%% >  ORDS1:=cuspORDSnotoo(epf1,C20,20);
%%           ORDS1 := [[0, 0],\quad, [1/2, -2],\quad, [1/4, -2],\quad, [1/5, 0],\quad, [1/10, 2],\quad],\quad
%% 
%% >  ORDS2:=cuspORDSnotoo(epf2,C20,20);
%%           ORDS2 := [[0, 1],\quad, [1/2, 0],\quad, [1/4, -1],\quad, [1/5, 0],\quad, [1/10, 1],\quad],\quad
%% 
%% > quit
%% memory used=2.7MB, alloc=8.3MB, time=0.11
% 
% end prog21.tex
%%%%%%%%%%%%%%%%%%%%%%%%%%%%%%%%%%%%%%%%%%%%%%%%%%%%%%%%%

        \vskip 10pt\noindent
{\it\footnotesize STEP 3b}. \quad  Use Theorem \thm{UpLB}
to find a lower bound for
$$
\ORD(\Up_{p} g, r, \Gamma_0(20)) 
$$
for each cusp $r$ of $\Gamma_0(20)$.
%%%%%%%%%%%%%%%%%%%%%%%%%%%%%%%%%%%%%%%%%%%%%%%%%%%%%%%%%
\showtfile{prog22.tex}
\mapleshade{
\sn
{\tt
\bmaplesess
\bmaple
$>$ \enskip  C20:=cuspmake(20) minus\{1/20\}: \newline
\emaple
\bmaple
$>$ \enskip  ORDSg:=[seq([cusp,UpLB(epg,cusp,100,5)], cusp in C20)],\quad;
\emaple
\bmapleout
$$
[[0,0],\quad[{1}/{2},-2],\quad[{1}/{4},-2],\quad[{1}/{5},-{1}/{5}],\quad[1/10,{3}/{5}] 
$$
\emapleout
\emaplesess
}\vskip0pt\noindent
}
%%%%%%%%%%%%%%%%%%%%%%%%%%%%%%%%%%%%%%%%%%%%%%%%%%%%%%%%%
%%     |\^/|     Maple 2018 (X86 64 WINDOWS)
%% ._|\|   |/|_. Copyright (c) Maplesoft, a division of Waterloo Maple Inc. 2018
%%  \  MAPLE  /  All rights reserved. Maple is a trademark of
%%  <____ ____>  Waterloo Maple Inc.
%%       |       Type ? for help.
%% >  with(qseries):
%% >  with(ETA):
%% >  gpg:=[100, -3, 50, 5, 25, -2, 10, -8, 5, 4, 4, 3, 2, 3, 1, -2],\quad:
%% >  epg:=gp2etaprod(gpg);
%%                               5           4           3           3
%%                    eta(50 tau)  eta(5 tau)  eta(4 tau)  eta(2 tau)
%%             epg := -------------------------------------------------
%%                                3            2            8         2
%%                    eta(100 tau)  eta(25 tau)  eta(10 tau)  eta(tau)
%% 
%% >  C20:=cuspmake(20) minus{1/20}:
%% >  ORDSg:=[seq([cusp,UpLB(epg,cusp,100,5)],\quad,cusp in C20)],\quad;
%%        ORDSg := [[0, 0],\quad, [1/2, -2],\quad, [1/4, -2],\quad, [1/5, -1/5],\quad, [1/10, 3/5],\quad],\quad
%% 
%% > quit
%% memory used=1.9MB, alloc=8.3MB, time=0.11
% 
% end prog22.tex
%%%%%%%%%%%%%%%%%%%%%%%%%%%%%%%%%%%%%%%%%%%%%%%%%%%%%%%%%

        \vskip 10pt\noindent
{\it\footnotesize STEP 4}. \quad  Calculate the constant $B$ in \eqn{UpB}.
%%%%%%%%%%%%%%%%%%%%%%%%%%%%%%%%%%%%%%%%%%%%%%%%%%%%%%%%%
\showtfile{prog23.tex}
\mapleshade
\sn
{\tt
\bmaplesess
\bmaple
$>$ \enskip  mintotGAMMA0ORDS([ORDSg,ORDS1,ORDS2],\quad,3);
\emaple
\bmapleout
$$
-{{18}/{5}} 
$$
\emapleout
\emaplesess
}\vskip0pt\noindent
%%%%%%%%%%%%%%%%%%%%%%%%%%%%%%%%%%%%%%%%%%%%%%%%%%%%%%%%%
%%     |\^/|     Maple 2018 (X86 64 WINDOWS)
%% ._|\|   |/|_. Copyright (c) Maplesoft, a division of Waterloo Maple Inc. 2018
%%  \  MAPLE  /  All rights reserved. Maple is a trademark of
%%  <____ ____>  Waterloo Maple Inc.
%%       |       Type ? for help.
%% >  with(qseries):
%% >  with(ETA):
%% >  gpf1:=[10, 8, 5, -4, 2, -8, 1, 4],\quad:
%% >  epf1:=gp2etaprod(gpf1):
%% >  gpf2:=[20, -3, 10, 5, 5, -2, 4, -1, 2, -1, 1, 2],\quad:
%% >  epf2:=gp2etaprod(gpf2):
%% >  C20:=cuspmake(20) minus {1/20}:
%% >  ORDS1:=cuspORDSnotoo(epf1,C20,20);
%%           ORDS1 := [[0, 0],\quad, [1/2, -2],\quad, [1/4, -2],\quad, [1/5, 0],\quad, [1/10, 2],\quad],\quad
%% 
%% >  ORDS2:=cuspORDSnotoo(epf2,C20,20);
%%           ORDS2 := [[0, 1],\quad, [1/2, 0],\quad, [1/4, -1],\quad, [1/5, 0],\quad, [1/10, 1],\quad],\quad
%% 
%% >  gpg:=[100, -3, 50, 5, 25, -2, 10, -8, 5, 4, 4, 3, 2, 3, 1, -2],\quad:
%% >  epg:=gp2etaprod(gpg):
%% >  ORDSg:=[seq([cusp,UpLB(epg,cusp,100,5)],\quad,cusp in C20)],\quad;
%%        ORDSg := [[0, 0],\quad, [1/2, -2],\quad, [1/4, -2],\quad, [1/5, -1/5],\quad, [1/10, 3/5],\quad],\quad
%% 
%% >  mintotGAMMA0ORDS([ORDSg,ORDS1,ORDS2],\quad,3);
%%                                      -18/5
%% 
%% > quit
%% memory used=2.0MB, alloc=8.3MB, time=0.08
% 
% end prog23.tex
%%%%%%%%%%%%%%%%%%%%%%%%%%%%%%%%%%%%%%%%%%%%%%%%%%%%%%%%%

\vskip 10pt\noindent
{\it\footnotesize STEP 5}. \quad  To prove the identity \eqn{Up5id} 
we need to
verify that 
$$
\ORD(h(\tau),\infty,\Gamma_0(20)) > 3,
$$
where
$$
h(\tau) = 
\Up_5(g(\tau)) - \left(5\,f_1(\tau) + 2\,f_2(\tau)\right).
$$
%%%%%%%%%%%%%%%%%%%%%%%%%%%%%%%%%%%%%%%%%%%%%%%%%%%%%%%%%
\showtfile{prog24.tex}
\mapleshade{
\sn
{\tt
\bmaplesess
\bmaple
$>$ \enskip  U5g:=sift(etaprodtoqseries(epg,1010),q,5,0,1000): \newline
\emaple
\bmaple
$>$ \enskip  h:=U5g - qetacombo(5*epf1 + 2*epf2,210): \newline
\emaple
\bmaple
$>$ \enskip  series(h,q,201);
\emaple
\bmapleout
$$
(O \left( {q}^{201} \right) ) 
$$
\emapleout
\emaplesess
}\vskip0pt\noindent
}

This completes the proof of the identity \eqn{Up5id}.            
We only had to show that the coefficient of 
$q^j$ was zero in the $q$-expansion of 
$h(\tau)$ for $j \le 4$.  We actually did it for 
$j \le 200$ as a check.

STEPS 1--5 may be automated using the following function.

\medskip

\texttt{provemodfuncGAMMA0UpETAid(EP,p,etacombo,N)} ---
%%\texttt{provemodfuncGAMMA0id(etaid,N)} --- 
attempts to prove the identity
$$
\Up_p(\mbox{\texttt{EP}}) = \mbox{\texttt{etacombo}},
$$
where \texttt{EP} is an eta-product and a modular function on $\Gamma_0(pN)$,
$p$ i sprime, $p\mid N$, and \texttt{etacombo} is a linear combination
of eta-products which are all modular functions on $\Gamma_0(N)$.
It returns the constant $B$ in equation \eqn{UpB}
and prints (if possible) the details of the verification and proof of the 
identity.
If \texttt{xprint=true} then more details of the
verification are printed. When this function is called there is a query asking
whether to verify the identity. Enter \texttt{yes} to carry out the verification.
%%%%%%%%%%%%%%%%%%%%%%%%%%%%%%%%%%%%%%%%%%%%%%%%%%%%%%%%%
\showtfile{prog25.tex}
\mapleshade{
\sn
{\tt
\bmaplesess
\bmaple
$>$ \enskip  etacombo:=5*epf1 + 2*epf2: \newline
\emaple
\bmaple
$>$ \enskip  provemodfuncGAMMA0UpETAid(epg,5,etacombo,20); \newline
\emaple
\bmapleout
*** There were NO errors. \newline
*** o EP is an MF on Gamma[0](100)\newline
*** o Each term in the etacombo is a  modular function on\newline
      Gamma0(20). \newline
*** o We also checked that the total order of\newline
      each term etacombo was zero.\newline
*** To prove the identity U[5](EP)=etacombo we need to show\newline
    that v[oo](ID) > 3    This means checking up to q\carrot(4).\newline
Do you want to prove the identity? (yes/no)\newline
You entered yes.\newline
We verify the identity to O(q\carrot(43)).\newline
We find that LHS - RHS is \newline
$$
O(q^{43})
$$
RESULT: The identity holds to O(q\carrot(43)).\newline
CONCLUSION: This proves the identity since we had only\newline
            to show that v[oo](ID) > 3.
\emapleout
\emaplesess
}\vskip0pt\noindent
}

\texttt{provemodfuncGAMMA0UpETAidBATCH(EP,p,etacombo,N)} --- 
is a version of\newline 
\texttt{provemodfuncGAMMA0UpETAid} that prints less detail and does not query.

%%%%%%%%%%%%%%%%%%%%%%%%%%%%%%%%%%%%%%%%%%%%%%%%%%%%%%%%%
\showtfile{prog27.tex}
\mapleshade{
\sn
{\tt
\bmaplesess
\bmaple
$>$ \enskip  provemodfuncGAMMA0UpETAidBATCH(epg,5,etacombo,20);\newline
\emaple
\bmapleout
 *** There were NO errors. \newline
 *** o EP is an MF on Gamma[0](100)\newline
 *** o Each term in the etacombo is a  modular function on\newline
       Gamma0(20). \newline
 *** o We also checked that the total order of\newline
       each term etacombo was zero.\newline
 *** To prove the identity U[5](EP)=etacombo we need to show\newline
     that v[oo](ID) > 3    This means checking up to q\carrot(4).\newline
 We find that LHS - RHS is \newline
$$
O \left( {q}^{43} \right) 
$$
$$
[1,-3,(O \left( {q}^{43} \right) )] 
$$
\emapleout
\emaplesess
}\vskip0pt\noindent
}

Let $L=[1,-3,O\left(q^{43}\right)]$. $L[1]=1$ means the identity is proved. $L[2]=B$ (the constant in equation \eqn{Bdef}), and  we see that $B=-2$. 
This confirms an earlier calculation using 
\texttt{mintotORDS}. $L[3]=O\left(q^{43}\right)$ mean the identity
was checked up to $q^{42}$.
To print out minimal information set \texttt{noprint:=true}.

%%%%%%%%%%%%%%%%%%%%%%%%%%%%%%%%%%%%%%%%%%%%%%%%%%%%%%%%%
\showtfile{prog28.tex}
\mapleshade{
\sn
{\tt
\bmaplesess
\bmaple
$>$ \enskip  noprint:=true: \newline
\emaple
\bmaple
$>$ \enskip  provemodfuncGAMMA0UpETAidBATCH(epg,5,etacombo,20);
\emaple
\bmapleout
$$
[1,-3,(O \left( {q}^{43} \right) )] 
$$
\emapleout
\emaplesess
}\vskip0pt\noindent
}


\begin{thebibliography}{10}
\bibitem{Bi89}
Anthony J.~F. Biagioli, \emph{A proof of some identities of {R}amanujan using
  modular forms}, Glasgow Math. J. \textbf{31} (1989), no.~3, 271--295.
  \MR{1021804}
\bibitem{Ch-La04}
Kok~Seng Chua and Mong~Lung Lang, \emph{Congruence subgroups associated to the
  monster}, Experiment. Math. \textbf{13} (2004), no.~3, 343--360. \MR{2103332}
\bibitem{Ga99b}
Frank Garvan, \emph{A {$q$}-product tutorial for a {$q$}-series {MAPLE}
  package}, S\'{e}m. Lothar. Combin. \textbf{42} (1999), Art. B42d, 27, The
  Andrews Festschrift (Maratea, 1998). \MR{1701583}
\bibitem{GoHu81}
B.~Gordon and K.~Hughes, \emph{Ramanujan congruences for {$q(n)$}}, Analytic
  number theory ({P}hiladelphia, {P}a., 1980), Lecture Notes in Math., vol.
  899, Springer, Berlin-New York, 1981, pp.~333--359. \MR{654539}
\bibitem{Li75}
G\'{e}rard Ligozat, \emph{Courbes modulaires de genre {$1$}}, Soci\'{e}t\'{e}
  Math\'{e}matique de France, Paris, 1975, Bull. Soc. Math. France, M\'{e}m.
  43, Suppl\'{e}ment au Bull. Soc. Math. France Tome 103, no. 3. \MR{0417060}
\bibitem{Ne59}
Morris Newman, \emph{Construction and application of a class of modular
  functions. {II}}, Proc. London Math. Soc. (3) \textbf{9} (1959), 373--387.
  \MR{0107629}
\bibitem{Onbook04}
Ken Ono, \emph{The web of modularity: arithmetic of the coefficients of modular
  forms and {$q$}-series}, CBMS Regional Conference Series in Mathematics, vol.
  102, Published for the Conference Board of the Mathematical Sciences,
  Washington, DC; by the American Mathematical Society, Providence, RI, 2004.
  \MR{2020489}
\bibitem{Ra}
Robert~A. Rankin, \emph{Modular forms and functions}, Cambridge University
  Press, Cambridge, 1977. \MR{0498390}
\end{thebibliography}
\end{document}